\documentclass[a4paper,10pt]{article}

\usepackage{amsbsy,amssymb,amsmath}
\usepackage{amsthm}  
\usepackage{caption,color,graphicx,graphics}
\usepackage{subfig}
\usepackage[a4paper,text={6.5in,9.5in},centering]{geometry}
\usepackage[colorlinks=true,hypertexnames=false,linkcolor=blue,citecolor=blue]{hyperref}
\usepackage[numbers,comma,square,sort&compress]{natbib}

\graphicspath{{figures/}}
\usepackage{psfrag}

\newtheorem{Lemma}{Lemma}[section]

\newtheorem{Hypothesis}{Hypothesis}[section]

\newtheorem{Theorem}{Theorem}[section]

\newcounter{hyp}
\newtheorem{hypothesis}{}

\newcounter{listitem}
\newenvironment{romanlist}{
\begin{list}{\hfill(\roman{listitem})}{\usecounter{listitem}
                 \itemindent0pt
                 \labelsep1.5ex
                 \leftmargin5.5ex
                 \labelwidth4ex
                 \topsep1ex
                 \itemsep0.5ex
                 \parsep0ex}}{%
\end{list}}

\makeatletter\@addtoreset{equation}{section}\makeatother

\newcommand{\R}{{\mathbb R}}

\newcommand{\N}{{\mathbb N}}

\newtheorem{Proof}{Proof}

\newenvironment{prooof}[1][.]%
 {\begin{trivlist}\item[]\textbf{Proof#1 }}%
 {\hspace*{\fill}$\blacksquare$\end{trivlist}}

\title{Nonreversible Homoclinic Snaking}

\author{J\"urgen Knobloch\\{\small  Institute of Mathematics, Ilmenau University of Technology} \\ {\small
e-mail: juergen.knobloch@tu-ilmenau.de}
\and
Thorsten Rie{\ss}\\{\small INCIDE, University of Konstanz}\\ {\small e-mail: thorsten.riess@uni-konstanz.de}
\and
Martin Vielitz\\{\small  Institute of Mathematics, Ilmenau University of Technology} \\ {\small
e-mail: martin.vielitz@tu-ilmenau.de}
}

\begin{document}
\date{\today}
\maketitle

\begin{abstract}
Homoclinic snaking refers to the sinusoidal ``snaking'' continuation curve of homoclinic orbits near a 
heteroclinic cycle connecting an equilibrium $E$ and a periodic orbit $P$. Along this curve the homoclinic orbit
performs more and more windings about the periodic orbit. Typically this behaviour appears in reversible Hamiltonian 
systems. Here we discuss this phenomenon in systems without any particular structure. We give a rigorous analytical verification of homoclinic snaking under certain assumptions on the behaviour of the stable and unstable manifolds of $E$ and $P$. We show how the snaking behaviour depends on the signs of the Floquet multipliers of $P$. Further we present a nonsnaking scenario. Finally we show numerically that these assumptions are fulfilled in a model equation.
\end{abstract}


{\bf AMS subject classifications.} global bifurcation, homoclinic snaking, heteroclinic cycle

\section{Introduction}

In this paper we study analytically a certain continuation scenario, the
so-called  {\sl Homoclinic Snaking}, of homoclinic orbits in systems without
particular structure such as reversibility or Hamiltonian structure.  Consider
an ordinary differential equation $\dot{x} = f(x)$ with $x\in{\mathbb R}^n$.
Given an equilibrium $E$, a homoclinic orbit to $E$ is a solution that converges
to $E$ as $t \to\pm \infty$. If $E$ is hyperbolic, a
homoclinic orbit lies in the intersection of the stable and unstable manifolds
$W^s(E)$ and $W^u(E)$ of the equilibrium.  For general differential equations,
stable and unstable manifolds of hyperbolic equilibria  will typically not
intersect  by the Kupka-Smale theorem \cite{palmel82}. Transversality arguments
show that in one-parameter families of differential equations one can expect an
intersection, and hence a homoclinic orbit, to occur persistently at an isolated
parameter value.  Therefore a continuation of a homoclinic orbit can be carried
out in two-parameter families of differential equations. 

Indeed, the homoclinic orbits we consider live within a small neighbourhood of a
heteroclinic cycle connecting a hyperbolic equilibrium $E$ and a hyperbolic
periodic orbit $P$. Such a cycle consists, besides $E$ and $P$, of two orbits
$\gamma_{\scriptscriptstyle{\rm EtoP}}$ and $\gamma_{\scriptscriptstyle{\rm
PtoE}}$ with $\lim_{t\to -\infty}\gamma_{\scriptscriptstyle{\rm EtoP}}(t)=E$,
$\lim_{t\to \infty}\gamma_{\scriptscriptstyle{\rm EtoP}}(t)=P$, and $\lim_{t\to
-\infty}\gamma_{\scriptscriptstyle{\rm PtoE}}(t)=P$, $\lim_{t\to
\infty}\gamma_{\scriptscriptstyle{\rm PtoE}}(t)=E$, respectively. The orbit
$\gamma_{\scriptscriptstyle{\rm EtoP}}$ is called heteroclinic orbit connecting
$E$ to $P$ or in short, EtoP connecting orbit , or just EtoP connection. A similar terminology holds for
$\gamma_{\scriptscriptstyle{\rm PtoE}}$ just with interchanging $E$ and $P$. The
complete cycle we call EtoP cycle.  More precisely, the homoclinic
orbits to the
equilibrium under consideration are one-homoclinic
orbits with respect to the
given EtoP cycle. That means they move once along the cycle before returning to
the equilibrium.

Replace the periodic orbit $P$ by a hyperbolic equilibrium $\hat E$, and assume that for a critical value in a
two-dimensional parameter space there exists a heteroclinic Eto$\hat{\rm E}$ cycle. Then, 
typically there is a curve emanating from the
critical value such that for all parameter values on this curve
there exits a homoclinic orbit to $E$ \cite{SSTC01,ChDeTe90}.
The homoclinic orbit spends more and more time near $\hat E$
when moving along this curve towards the critical point.

Considering homoclinic orbits to $E$ in a neighbourhood of an EtoP
cycle, we also find that along the continuation curve the
homoclinic orbits spend more and more time near the periodic
orbit. In contrast to the Eto$\hat{\rm E}$-case
the continuation curve does not converge to a point, but it
approaches a curve segment as $h_1^b$ does
in Figure~\ref{f:laser} below. This makes the consideration global in the
parameter space. Before discussing this in more detail we
consider the problem from the homoclinic snaking point of view.

\begin{figure}[ht]
\centering\input{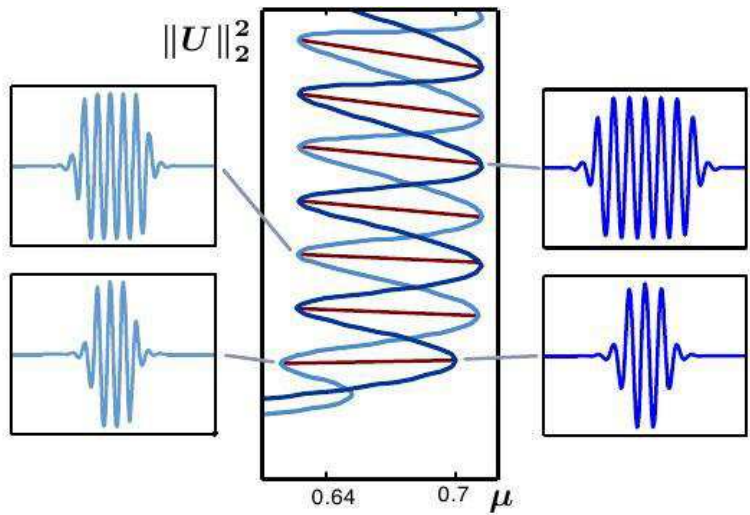}
\caption{Bifurcation diagram for homoclinic orbits of the steady states $-(1+\partial_x^2)^2U-\mu U+2 U^3-U^5=0$ of the Swift-Hohenberg equation. The central panel shows the typical snakes and ladder structure. The outer panels show the $U$ vs. time plot for the corresponding points on the snaking curves.}\label{f:sw_ho}
\end{figure}

In the context of ordinary differential equations the notion {\sl Homoclinic Snaking} originally denotes a  continuation scenario of homoclinic orbits in reversible Hamiltonian systems. In Hamiltonian systems the situation is somewhat different to the one described above. Both the stable and the unstable manifold of a hyperbolic equilibrium are in the same levelset of the Hamiltonian. Therefore they will typically intersect transversely (within this levelset) and a homoclinic orbit can be expected to occur persistently in a single system. Hence a continuation can be done in one-parameter families.
Typical continuation curves related to a homoclinic snaking scenario are displayed in Figure~\ref{f:sw_ho}, cf. \cite[Figure~1.1]{BKLSW:08}. 

Figure~\ref{f:sw_ho} shows continuation curves of homoclinic orbits related to the fourth order scalar equation
\begin{equation}\label{e:sw_ho}
 -(1+\partial_x^2)^2U-\mu U+2 U^3-U^5=0.
\end{equation}
Considered as a system in $\R^4$, equation  \eqref{e:sw_ho} is
a reversible Hamiltonian system with family parameter $\mu$.  The two
intertwined wiggly curves correspond to symmetric (w.r.t. the reflection
$x\mapsto -x$ and $U\mapsto -U$) homoclinic orbits asymptotic to a hyperbolic
equilibrium. These curves are also called {\sl snaking curves}.  The term
snaking is due to the sinusoidal shape of the continuation curves. Indeed, the
homoclinic orbits under consideration are one-homoclinic w.r.t. a symmetric
heteroclinic cycle connecting a symmetric equilibrium with a symmetric periodic
orbit. Restricted to the corresponding levelset of the
Hamiltonian, the periodic orbit is hyperbolic.
The excursion of the homoclinic orbit to the periodic orbit
lasts longer and longer along the continuation curves.
The homoclinic orbit performs more and more windings about the
periodic orbit -- this corresponds to the
increase of its $L^2$-norm.  Roughly
speaking, the $\mu$-range of the snaking curves is the
$\mu$-range for which the heteroclinic cycle does exist.  At the endpoints of
the $\mu$-interval the involved EtoP and PtoE connections (note that they are
images of each other by the reversing symmetry) simultaneously undergo a
saddle-node bifurcation.  Note that the snaking curves indicate saddle-node
bifurcations of the symmetric homoclinic orbits. These bifurcation points are
close to the endpoints of the $\mu$-range for the heteroclinic cycle. In addition to the
snaking curves, Figure~\ref{f:sw_ho} also displays a ladder
structure. The rungs connecting the two snaking curves
correspond to asymmetric homoclinic orbits to the equilibrium. The asymmetric
homoclinic orbits bifurcate from the symmetric ones via pitchfork bifurcation.
These bifurcation points are close to the saddle-nodes of the symmetric
homoclinic orbits.  However, in the context of the present treatment it is
enough to focus on the features of one single snaking curve.  Bifurcation
diagrams as displayed in Figure~\ref{f:sw_ho} have been discussed for instance
in \cite{BKLSW:08,WoCh:99,CoRiTr:00,BuKn:06}. For a more complete list of
references we refer to \cite{BKLSW:08},
but with the addition that
homoclinic snaking, also called collapsed snaking, near Eto$\hat{\rm E}$~cycles in one-parameter families of reversible systems has been studied in \cite{KnoWag05,KnoWag08}.

More recently, Krauskopf, Oldeman and Rie{\ss} \cite{KrauOld06,KrauRie:08}
numerically discovered a similar effect in a
system without any particular structure such as
reversibility or Hamiltonian structure.  The corresponding family of vector fields in $\R^3$ with family parameter $(\nu_1,\nu_2)$ can be written in the form
\begin{equation}\label{e:laser}
\left.
 \begin{array}{ll}
 \dot{x} &\!\!\!\!\!= \nu_1 x -y + x  \sin \varphi - (x^2 + y^2)x+ 0.01 ( 2\cos\varphi + \nu_2)^2 
\\[1ex]
\dot{y} &\!\!\!\!\!= \nu_1 y +  x + y  \sin \varphi - (x^2 + y^2)y +
0.01\pi ( 2\cos\varphi +\nu_2)^2 
\\[1ex] 
\dot{\varphi} &\!\!\!\!\!= \nu_2 -(x^2+y^2)+2\cos\varphi
\end{array}
\right\} =:F(x,y,\varphi,\nu_1,\nu_2).
\end{equation}

Figure~\ref{f:laser} displays a continuation curve $h_1^b$ for a homoclinic
orbit detected numerically in \cite{KrauOld06,KrauRie:08}.  Again, and not only
due to its shape, we address this curve as snaking curve.
\begin{figure}[hb]
\centering
\includegraphics[scale=0.8]{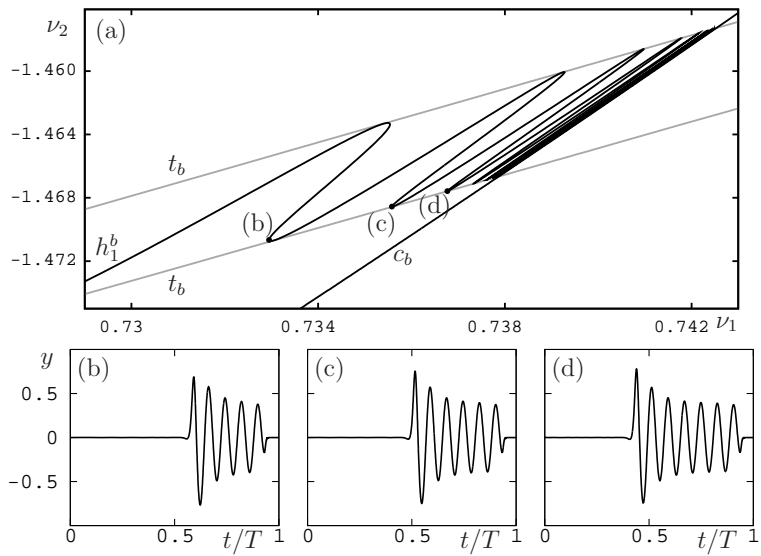} 
\caption{Snaking diagram of a three-dimensional laser model \eqref{e:laser}. Panel (a) shows the snaking curve $h_1^b$ together with the continuation curve $c_b$ of the PtoE connecting orbits and the locus $t_b$ of saddle-nodes of the EtoP connecting orbits. Panels (b) - (d) show $y$ vs. time plots of the $y$ component  at the corresponding points (b) - (d) in panel (a).}\label{f:laser}
\end{figure}
As in the reversible case, the homoclinic orbit under consideration is one-homoclinic w.r.t. a
heteroclinic cycle connecting a hyperbolic equilibrium with a hyperbolic
periodic orbit, and along the continuation curve the
homoclinic orbits performs more and more windings about the periodic orbit
-- cf. panels (b)-(d) in
Figure~\ref{f:laser}, which show plots of one state variable
corresponding to points indicated in panel (a).  The plot of the $L^2$-norm
of the $(x,y)$-part of the solution versus the parameter
$\nu_1$ or $\nu_2$, respectively, behaves as in the
Hamiltonian case, cf.  Figure~\ref{f:laser1}. The snaking
behaviour w.r.t. both parameters is due to the declination of the
curve $c_b$ in Figure~\ref{f:laser}\,(a). This defines
intervals within which the parameters move while $h_1^b$ approaches $c_b$.
We refer also to Section~\ref{s:numver} for more numerical
results regarding this system.

\begin{figure}[ht]
\centering
\includegraphics[scale=0.6]{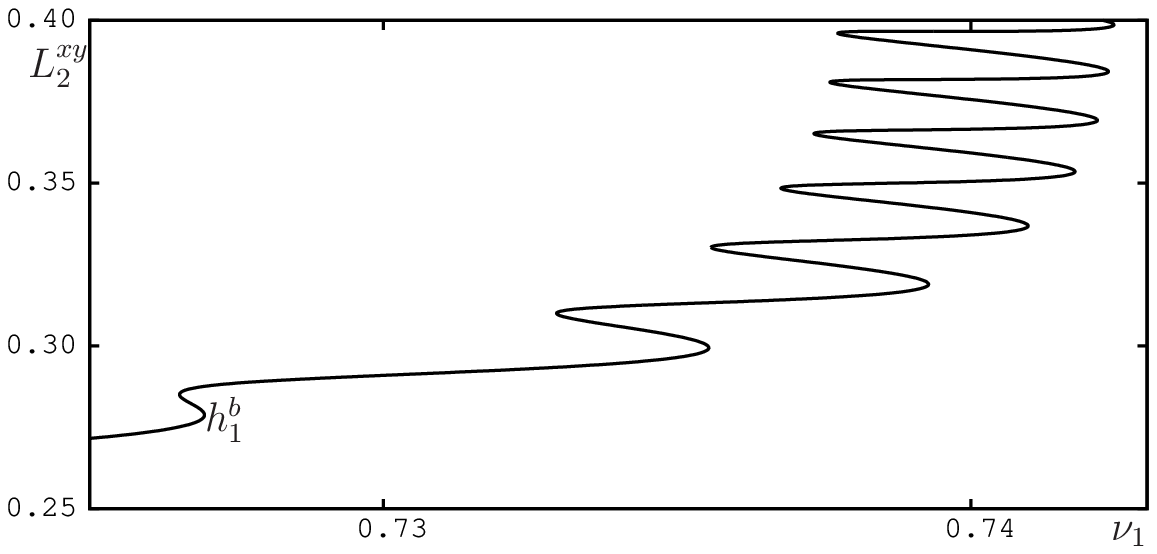} \quad\quad
\includegraphics[scale=0.6]{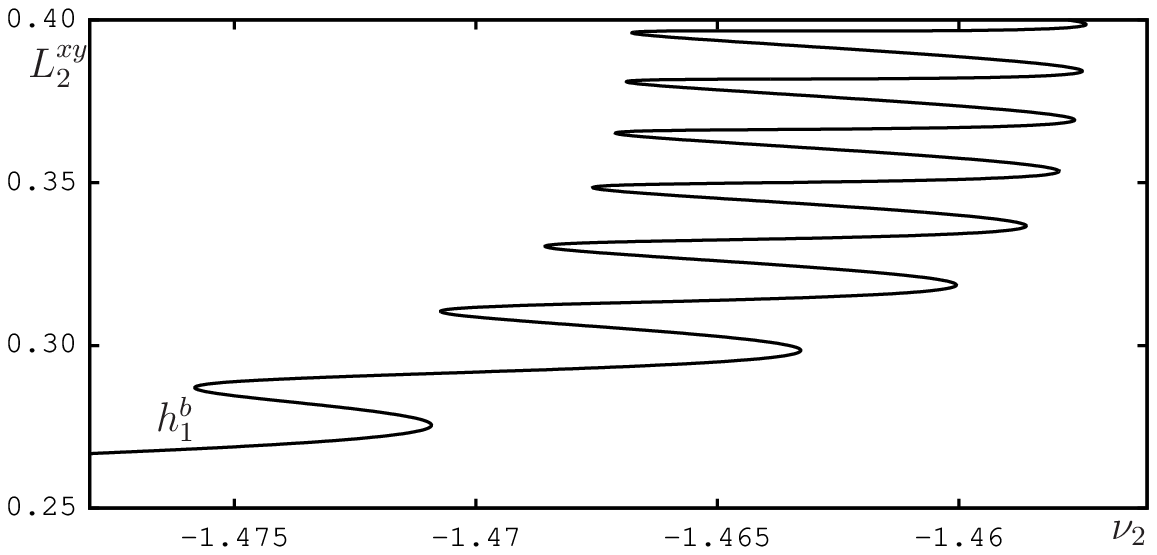} 
\caption{Snaking diagrams of a three-dimensional model \eqref{e:laser}. Shown are plots of the $L^2$-norm of $(x,y)$ vs. $\nu_1$ and $\nu_2$, respectively, along the snaking curve $h_1^b$.
}\label{f:laser1}
\end{figure}

Although system \eqref{e:laser} has the same dimension as the restriction of the
above Hamiltonian system to a levelset, the geometry is quite different.  Here
one of the heteroclinic connections constituting the EtoP cycle does not lie in
a transversal intersection of the corresponding stable and unstable manifolds.
Hence, by the same arguments as given above for homoclinic
orbits, one can expect to find it on a continuation curve in the
$(\nu_1,\nu_2)$-space. This is the curve $c_b$ in Figure~\ref{f:laser}\,(a). The
other connecting orbit is robust -- as the ones in the
Hamiltonian case. This connecting orbit exists within the stripe delimited by
the curves $t_b$, cf. again Figure~\ref{f:laser}\,(a).  In other words, the
region where the snaking curve is located is related to the existence of the
robust heteroclinic connection between the equilibrium and the periodic orbit.
Remarkably, the snaking curve accumulates at the
curve segment defined by the intersection of the curve $c_b$
with the stripe delimited by the curves $t_b$.  In other words, the snaking
curve accumulates at the line segment for which the EtoP cycle does exist. 

In this paper we give an analytical verification of the bifurcation diagram in
Figure~\ref{f:laser} within a more general setup.  
We consider a general two-parameter family of differential equations in
$\R^3$. In this general context we denote the family parameter by
$(\lambda_1,\lambda_2)$ taken from a closed rectangle $Q$.  We assume that there
is a closed interval $I_2$ such that for all $(\lambda_1,\lambda_2)\in
\{0\}\times I_2$ there exists an EtoP cycle connecting an equilibrium and a
periodic orbit, both are assumed to be hyperbolic, cf.
Hypothesis~\ref{h:general}. Let $W^s(E)$ be two-dimensional.  
We want to note that this setting is related to the vector field $-F$ in \eqref{e:laser}. Compare also the caption of Figure~\ref{f:laser}.
Further we
constitute conditions concerning the global intersection behaviour of $W^u(P)$,
$W^s(E)$ and $W^u(E)$, $W^s(P)$, respectively, cf.  Hypotheses~\ref{h:PtoE1} --
\ref{h:EtoP}. Note that all manifolds depend on $\lambda$,
which is so far suppressed from the notation.  Indeed, the snaking behaviour is
mainly influenced by the behaviour of the intersection of $W^u(P)$ and $W^s(E)$.
Consider a small torus ${\mathcal T}$ around $P$, and let $\Sigma^{out}$ be a
small stripe on this torus around $W^u(P)\cap {\mathcal T}$.  Similarly we define $\Sigma^{in }$
as a small stripe on this torus around $W^s(P)$.  In $\Gamma\subset S^1\times Q$
we collect the intersections of $W^u(P)$ and $W^s(E)$ in $\Sigma^{out}$
depending on $\lambda$: 
\[ 
\Gamma:=\{(\varphi,\lambda): W^s(E,\lambda)\cap
W^{uu}(P(\varphi,\lambda),\lambda)\cap \Sigma^{out}\not=\emptyset\}, 
\] 
where
$W^{uu}(P(\varphi,\lambda))$ is the strong unstable fibre of
$P(\varphi,\lambda)\in P$.  The assumption that $\Gamma$ is
graph of a function $z=z(\varphi,\lambda_1)$ is essential for
the snaking behaviour.  Our hypotheses on $z$ imply amongst others that at the
endpoints of $I_2$ the PtoE connection undergoes a saddle-node bifurcation, and
what is more, for each $\lambda\in\{0\}\times {\rm int}\,I_2$
there are at least two EtoP cycles.

To get a better idea of what these assumptions include, suppose for simplicity
that the trace of $W^s(E)$ in $\Sigma^{out}$ is a closed curve which is simple
over $W^u(P)$.  So $W^s(E)$ can be seen as function $Z$ of values in $W^u(P)$.
To simplify matters further we assume that by changing $\lambda$ these curves
will be shifted against each other without changing their shape. Then $Z$ and
$z(\cdot,0)$ are directly related.

For $\lambda_1=0$ and all $\lambda_2\in I_2$, there is a unique
EtoP connection; all those $\lambda_2$ can be written in the
form $\lambda_2=z(\varphi,0)$.  Our hypotheses on the EtoP connection
provide the existence of a $2\pi$-periodic function
$\varphi_0^*(L)$ defining the base point of the strong stable fibre of $P$ which
intersects in $\Sigma^{in}$ the EtoP connection related to $\lambda_1=0$,
$\lambda_2=z(\varphi_0^*(L)+2L,0)$.

Our main snaking result, cf. Theorem~\ref{t:snaking}, says that all
one-homoclinic orbits near the primary EtoP cycles lie on one continuation
curve, which we refer to as the snaking curve. This curve can be parametrised by
the flight time $L$ of the one-homoclinic orbit between $\Sigma^{in}$ and
$\Sigma^{out}$.  For $L\to\infty$, this curve accumulates at
$\{0\}\times I_2$.  It turns out that the shape of the snaking curve is mainly
determined by $z$ -- more precisely $\lambda_2(L)$ arises as
a perturbation of ${\rm graph}\,z(\varphi_0^*(L)+2L,0)$.

In the scenario covered by Theorem~\ref{t:snaking} the nontrivial Floquet multipliers of $P$ are positive. 
Indeed, the snaking behaviour depends on the sign of these multipliers. If they are positive, the local (un)stable manifold of $P$ is topologically a cylinder, while for negative multipliers these local manifolds are topologically a M\"obius strip. Therefore, for positive multipliers both $\Sigma^{in}$ and $\Sigma^{out}$ consist of two connected components. In our analysis however, only one of these components, in each case, plays a role -- and there is only one way for the transition from $\Sigma^{in}$ to $\Sigma^{out}$.
But if the multipliers are negative, then $\Sigma^{in}$ and $\Sigma^{out}$ are connected, and both are winding twice around $P$ on ${\mathcal T}$. This  results in two different ways for the transition from $\Sigma^{in}$ to $\Sigma^{out}$, and this causes the existence of two snaking curves approaching $\{0\}\times I_2$ from different sides, cf. Theorem~\ref{t:snaking_floquet}.

Next we abandon our assumption on $\Gamma$ being graph of a function and assume
rather that for fixed $\lambda_1$ the set $\Gamma$ is a closed curve, and
replace Hypothesis~\ref{h:PtoE2} by Hypothesis~\ref{h:PtoE2-non-snake}.  This
prevents snaking. In this case no longer all homoclinic orbits are on one
continuation curve, instead there exists a sequence of closed
homoclinic continuation curves in the $\lambda$-plane accumulating at $\{0\}\times
I_2$, cf. Theorem~\ref{t:nonsnaking}.

Numerically we show that for the motivating system \eqref{e:laser} our
hypotheses generating snaking are fulfilled. In particular we verify
Hypothesis~\ref{h:PtoE2}.

For our analysis we use Fenichel coordinates near the periodic orbit, and within
this setting we consider solutions of a Shilnikov problem which we glue together
with the stable and unstable manifolds of the equilibrium. This procedure is the
same as the one utilised in \cite{BKLSW:08} and \cite{KLSW:10} to study
reversible Hamiltonian systems. To our knowledge there are only a few further
works presenting analytical results for the dynamics near EtoP cycles. In
\cite{KnoRie:10} a Lin's method approach has been extended to treat heteroclinic
chains involving period orbits. These results are applied to EtoP cycles, in
particular to detect nearby one-homoclinic orbits. However, the results are more
local in nature. These results concern, in the context of Figure~\ref{f:laser},
the existence of one-homoclinic orbits for parameter values in the neighbourhood
of certain points on the curve $c_b$ -- but not in the
neighbourhood of an entire segment of $c_b$ as in the present paper.  In
\cite{Rademacher2005,Rademacher2008} a somewhat different (in handling the flow
near $P$) Lin's method approach has been used to study EtoP cycles of
codimension-one and codimension-two. In this language the EtoP cycle considered
in the present paper are of codimension-one. Rademacher's results about
homoclinic orbits near codimension-one EtoP cycles are of the same nature as the
ones in \cite{KnoRie:10}.  

In our analysis we assume the existence of a primary
EtoP cycle, and we make assumptions about its global (in parameter space)
behaviour.  In the Swift-Hohenberg equation the existence of heteroclinic
connections has been investigated analytically in \cite{ChKo:09,KoCh:06}.  

In \cite{CKKOR09}, amongst others, homoclinic snaking caused by an EtoP cycle in systems in $\R^3$ is
considered.
Using a combination of geometric and analytical arguments, the
snaking behaviour as displayed in Figure~\ref{f:laser} is explained.  More
precisely, based on a leading term approximation of the bifurcation
equation, one-homoclinic orbits near the saddle-node points of
the EtoP cycles are determined. Then the entire snaking curve is deduced by
using geometric arguments.  This has been done for both, when $P$ has
positive or negative Floquet multipliers.  Here in the present paper on the
contrary, we give a rigorous analytical verification of these scenarios. 

Numerically the homoclinic snaking scenario in the addressed model, which is inspired by semiconductor laser dynamics, has
been considered in several papers.  The snaking curve $h_1^b$, cf.
Figure~\ref{f:laser}, was first revealed in \cite{KrauOld06}. In
\cite{KrauRie:08} the  relation to the organising EtoP cycle was numerically
discovered in a bifurcation diagram similar to that in Figure~\ref{f:laser}.
This system was further investigated in \cite{CKKOR09,KnoRie:10}.  In
\cite{CKKOS07,CKKOR09} a similar snaking behaviour was (numerically) observed in
a nine-dimensional model equation of intracellular calcium dynamics.  The
remarkable feature in the bifurcation diagram is that turning points of the
snaking curve accumulate on six different values (and not on only two as in the
one displayed in Figure~\ref{f:laser}).  Note that our analysis is carried out
only for the case of a three-dimensional state space.  However, the set $\Gamma$
remains a curve also in higher dimensions. Then a corresponding function $z$ can
be defined, and the addressed feature can be explained by the number of critical
points of $z$, cf. also Figure~\ref{f:z_lambda_relat} below.  

In \cite{BuHoKn:09} homoclinic snaking in the transition from reversible
Hamiltonian systems to general systems using the example of the Swift-Hohenberg
equation is considered. Indeed this happens in $\R^4$, but the observed snaking
or nonsnaking behaviour in the perturbed system, respectively, discloses
features we discuss for general systems in $\R^3$.
In \cite[Figure~6]{BuHoKn:09} snaking curves of two different homoclinic orbits
are shown -- each curve displaying a behaviour as shown in
Figure~\ref{f:laser1} of the present paper. These two different homoclinic
orbits are remains of the unperturbed reversible Hamiltonian system.  Isolas of
homoclinic orbits as shown in \cite[Figures~2~and~3]{BuHoKn:09} are discussed in
Section~\ref{s:nonsnak_ana} of this paper.  In symmetric systems those isolas
are also observed in \cite[Figure~24]{BuKn:06}.  However, it is not the aim of
this paper to explain those transition processes.

This paper is organised as follows. In Section~\ref{s:setup_res} we present our hypotheses and formulate the main snaking result, Theorem~\ref{t:snaking}.
The proof of Theorem~\ref{t:snaking}, is then carried out in Section~\ref{s:snak_ana}.
In Section~\ref{s:neg_Floqu_mul} we treat negative Floquet multipliers. The results are summarised in Theorem~\ref{t:snaking_floquet}.
Afterwards we discuss one possible nonsnaking scenario in Section~\ref{s:nonsnak_ana}; Theorem~\ref{t:nonsnaking} covers the results of this section. 
In Section~\ref{s:numver} we verify numerically Hypothesis~\ref{h:PtoE2}, the main snaking assumption, in the laser model \eqref{e:laser}.

\section{Setup and main results}\label{s:setup_res}
We consider a smooth family of differential equations
\begin{equation}\label{e:system}
 \dot x=f(x,\lambda),\quad x\in\R^3, \quad \lambda=(\lambda_1,\lambda_2)\in Q\subset\R^2,
\end{equation}
where $Q=J_1\times J_2$ is a closed rectangle with $0\in \,{\rm int}\,Q$;
$J_1$, $J_2$ are closed intervals.

We assume the following
\begin{Hypothesis}\label{h:general}
 \begin{romanlist}
 \item $f(0,\lambda)\equiv 0$; The equilibrium $E:= \{x=0\}$ is hyperbolic, and $\dim W^u(E,\lambda)=1$, $\dim W^s(E,\lambda)=2$.
\item For all $\lambda\in Q$ there is a hyperbolic periodic orbit $P$. Further let
$\dim W^u(P,\lambda)=2$, $\dim W^s(P,\lambda)=2$.
For all $\lambda$ the minimal period of $P$ is $2\pi$.
\item There is a maximal interval $I_2:=[\underline{\lambda}_{\,2},\overline{\lambda}_2]\subsetneq J_2$, $\underline{\lambda}_{\,2}<\overline{\lambda}_2$, such that 
for $\lambda\in\{0\}\times I_2$
there is a heteroclinic cycle connecting $E$ and $P$.
\end{romanlist}
\end{Hypothesis}
The constant minimal
period can always be achieved by an appropriate time transformation.  The
interval $I_2$ is maximal in the sense that for $(\lambda_1=0,\lambda_2)$ and
$\lambda_2>\overline{\lambda}_2$ or $\lambda_2<\underline{\lambda}_{\,2}$ there
in no complete cycle. More precisely with our choice of dimensions, typically
the EtoP connection is of codimension-one -- that means it
appears along a curve in parameter space. This curve
is the $\lambda_2$-axis and the connection
splits up when moving off the $\lambda_2$-axis.  On the other
hand, the PtoE connection is typically robust. Nevertheless,
by changing parameters within a wider range this connection can
disappear, for instance in the course of a saddle-node
bifurcation. These scenarios are made more precise
by additional hypotheses below.

The three-dimensional state space enforces that both nontrivial Floquet multipliers of $P$ have the same sign, cf. \cite{SSTC98}.
\begin{Hypothesis}\label{h:pos_Floquet_mul}
 The nontrivial Floquet multipliers of $P$ are positive.
\end{Hypothesis}

The positivity of the Floquet multipliers of $P$ is exploited in the Fenichel
normal form near $P$, cf. Lemma~\ref{l:fenichel} below.  However, in Section~\ref{s:neg_Floqu_mul} we relaxe this hypothesis.

The following lemma can be seen as a motivation for our further considerations.
Roughly speaking, it says that under certain transversality
conditions on each curve $\kappa$ intersecting $\{\lambda_1=0\}$
transversely, there is a sequence of parameter values
accumulating at $\{\lambda_1=0\}$ for which a one-homoclinic orbit to the
equilibrium does exist, cf. \cite[Corollary 4.3]{KnoRie:10}.

\begin{Lemma}[\cite{KnoRie:10}]\label{l:accum}
Assume Hypotheses~\ref{h:general} and \ref{h:pos_Floquet_mul}, and
let $\kappa=\kappa(\mu)$ be a smooth curve in $Q$ intersecting $\{\lambda_1=0\}$ transversely in $(0,\hat\lambda_2)$, where $\hat\lambda_2\in(\underline{\lambda}_2,\overline{\lambda}_2)$ and $\kappa(0)=(0,\hat\lambda_2)$.
Assume further
\begin{romanlist}
 \item $\mathop\bigcup\limits_\mu \big(W^u(E,\kappa(\mu))\times\{\mu\}\big)
\pitchfork 
\mathop\bigcup\limits_\mu \big(W^s(P,\kappa(\mu))\times\{\mu\}\big)$
\item $W^s(E,(0,\hat\lambda_2))\pitchfork W^u(P,(0,\hat\lambda_2))$
\end{romanlist}
Then there is a sequence $(\mu_n)$, $\lim\limits_{n\to\infty}\mu_n=0$ such that for all $\lambda=\kappa(\mu_n)$, $n\gg 1$, there is a one-homoclinic orbit to $E$.
\end{Lemma}

Assumption (i) of the lemma claims that the extended unstable and stable manifold of the equilibrium and of the periodic orbit, respectively, intersect transversely, while assumption (ii) claims that the stable and unstable manifold of the equilibrium and the periodic orbit, respectively, intersect transversely.

\begin{figure}[h]
\centering
\input{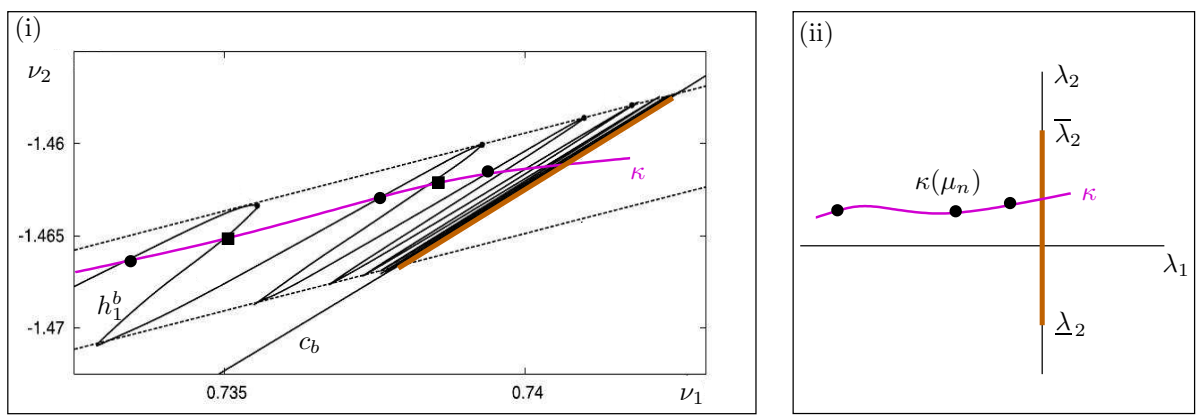}
\caption{One-homoclinic orbits on a curve $\kappa$ which intersects the continuation curve of the codimension-one heteroclinic orbits transversely. Panel~(i) is related to the laser model \eqref{e:laser}. The dots and squares indicate that the homoclinic orbits correspond to different EtoP cycles which exist at the intersection point of $\kappa$ and $c_b$. Panel~(ii) visualises the statement of Lemma~\ref{l:accum}.
}\label{f:accum}
\end{figure}

Now arises the question whether all $\kappa(\mu_n)$ lie on one
continuation curve as in our motivating example -- cf.
Figure~\ref{f:accum}. In panel (i) of this figure, the black dots and squares
correspond to parameter values on $\kappa$ for which a homoclinic orbit exists.
The different shapes indicate that the homoclinic orbits are related
to different EtoP cycles. Indeed, our
considerations in Section~\ref{s:numver} confirm that in the
intersection of $\kappa$ with $c_b$ there exist two EtoP cycles. This feature
has not been considered in panel (ii).

Denote the Floquet exponents of $P$ by
$\alpha^s(\lambda)<0<\alpha^u(\lambda)$ and assume the associated Floquet
multipliers to be positive. Moreover, let $\delta >0$ be a
sufficiently small constant and $I:=[-\delta,\delta]$, and let
$S^1:=\R/_{\sim_{2\pi}}$, and $x\sim_{2\pi} y \Leftrightarrow x=y\!\mod\!2\pi$.
We introduce the
so-called Fenichel coordinates $v=(v^c,v^s,v^u)\in S^1 \times I \times I$, which
are defined in a $\delta$-neighbourhood of the periodic orbit
$P$.  These coordinates are specially tailored to the hyperbolic structure of
$P$. This is reflected by the fact that the stable manifold $W^s(P,\lambda)$
in these coordinates corresponds to the set $\{v^u=0\}$,
whereas $\{v^s=0\}$ represents the unstable manifold $W^u(P,\lambda)$. Further
fixing $v^c=\varphi$ yields the single strong stable and strong unstable fibres
$W^{ss}(P(\varphi,\lambda),\lambda)$ and $W^{uu}(P(\varphi,\lambda),\lambda)$,
respectively. The periodic orbit itself is given by the set $\{v^s=0, \
v^u=0\}$.

\begin{Lemma}\label{l:fenichel}
Assuming Hypotheses \ref{h:general}\,(ii) and \ref{h:pos_Floquet_mul} are met, there is a smooth change of coordinates locally near $P$ such that $\dot{x}=f(x,\lambda)$ becomes
\begin{equation}\label{e:fenichel}
 \begin{array}{lll}
 \dot v^c&=& 1+A^c(v,\lambda)v^sv^u,
\\
 \dot v^s&=& (\alpha^s(\lambda)+A^s(v,\lambda))v^s,
\\
 \dot v^u&=& (\alpha^u(\lambda)+A^u(v,\lambda))v^u,
\end{array}
\end{equation}
where $v=(v^c,v^s,v^u)\in S^1 \times I \times I$ and $A^c, \ A^s, \ A^u$ are some smooth functions in $v$ and $\lambda$ with
\[
 A^i(v^c,0,0,\lambda)=0, \quad i=c,s,u, \quad \forall \, \lambda \in Q.
\]
\end{Lemma}

The Fenichel coordinates were introduced at first by Fenichel in the context of slow/fast systems, \cite{Fen79}. In \cite{Jon95} the transformation into the Fenichel coordinates in the context of slow/fast systems is described in more detail. 
The transformation near the hyperbolic periodic orbit $P$, and hence the proof of Lemma~\ref{l:fenichel}, is done in an analogous way.
The idea of the proof of Lemma \ref{l:fenichel} can also be found in \cite{BKLSW:08}.
However, note that in \cite{Fen79} and  \cite{Jon95} the Fenichel coordinates are derived merely locally, but for our purpose we need a global change of coordinates with respect to the periodic orbit $P$. To guarantee that the Fenichel coordinates can be introduced along the whole periodic orbit, we exploit the assumption that both nontrivial Floquet multipliers are positive,  Hypothesis~\ref{h:general}\,(ii), since this implies that the stable and unstable vector bundles of $P$ are orientable. 

\begin{figure}[ht]
\centering\input{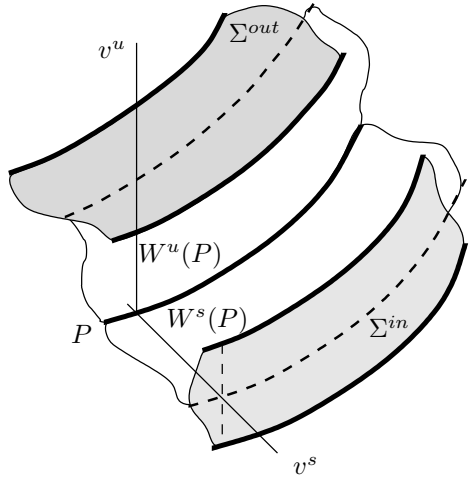}
\caption{The cross-sections $\Sigma^{in}$ and $\Sigma^{out}$.}\label{f:sigmas}
\end{figure}

Next we introduce sections near $P$.
\begin{equation}\label{e:def_Sigma}
 \Sigma^{in}:= S^1\times\{v^s=\delta\}\times I, \qquad
\Sigma^{out}:= S^1\times I\times\{v^u=\delta\},
\end{equation}
which are illustrated in Figure~\ref{f:sigmas}. Indeed, these sections are connected components of the sections $\Sigma^{in}$ and $\Sigma^{out}$ from the Introduction.
Further we define 
\begin{equation}\label{e:Gamma_def}
 \Gamma:=\{(\varphi,\lambda)\in S^1\times Q: W^s(E,\lambda)\cap W^{uu}(P(\varphi,\lambda),\lambda)\cap \Sigma^{out}\not=\emptyset\}.
\end{equation}
Thus $\Gamma$ consists of all the tuples $(\lambda, \varphi)$ for which there exists a PtoE connection that contains the strong unstable fibre $W^{uu}(P(\varphi,\lambda),\lambda)$ to the base point $P(\varphi,\lambda)$. 

Further, let $U_\Gamma$ be an open neighbourhood of $\Gamma$
in $S^1\times Q$.
\begin{Hypothesis}\label{h:PtoE1}
 There is a smooth function $g:U_\Gamma\to I$ and an $\epsilon>0$ such that
\begin{equation*}
  \{(\varphi,v^s,\delta)\in W^s(E,\lambda)\cap\Sigma^{out}: |v^s|<\epsilon, (\varphi,\lambda)\in U_\Gamma\}
=\{(\varphi,g(\varphi,\lambda),\delta): (\varphi,\lambda)\in U_\Gamma\}.
\end{equation*}
\end{Hypothesis}

\begin{figure}[hb]
\centering\input{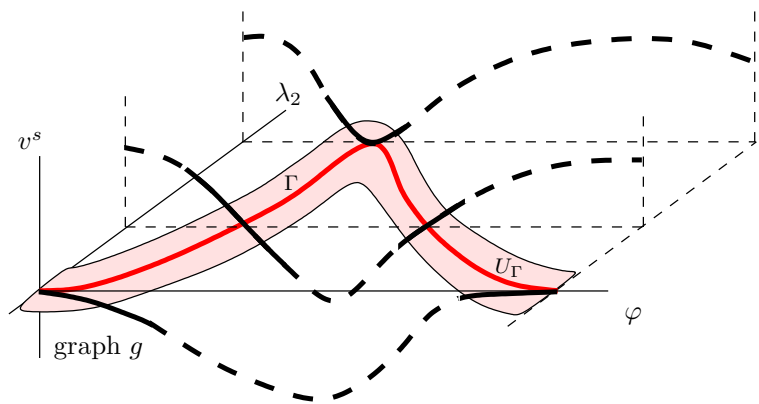}
 \caption{Visualisation of Hypothesis~\ref{h:PtoE1} and its consequence: In this illustration $\lambda_1$ is fixed with $\lambda_1=0$. $\Gamma$ coincides with the zeros of $g$. The graph of $g$ is only drawn for a sample of $\lambda_2$-values. The dashed lines indicate a possible continuation of $g$ outside of $U_\Gamma$.
}\label{f:graph_g}
\end{figure}

As a consequence of that hypothesis, we get that $\Gamma$ coincides with the zeros of $g$, cf. Figure~\ref{f:graph_g}:
\begin{equation}\label{e:Gamma_rep}
 \Gamma:=\{(\varphi,\lambda)\in S^1\times Q: g(\varphi,\lambda)=0\}.
\end{equation}

\begin{figure}[ht]
\centering\input{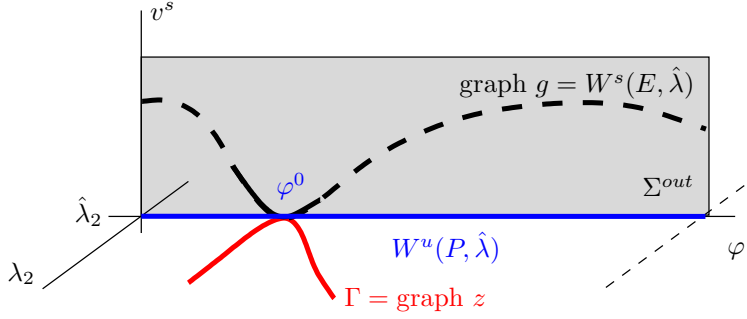}
 \caption{Quadratic tangency of $W^u(P)$ and $W^s(E)$ as a consequence of Hypothesis~\ref{h:PtoE3}: In this illustration $\lambda_1$ is fixed with $\lambda_1=\hat\lambda_1$.
 For $\hat\lambda_1=0$, this is an enlargement of a detail around the turning point of $\Gamma$ in Figure~\ref{f:graph_g}.
}\label{f:PtoE3_geom}
\end{figure}

Figure~\ref{f:graph_g} does already include some specific features of $g$ or $\Gamma$, respectively, which we demand in the following hypothesis:
\begin{Hypothesis}\label{h:PtoE2}
 \begin{romanlist}
 \item There is a constant $b>0$ such that $|g_{\lambda_2}(\varphi,\lambda)|\geq b$, for all $(\varphi,\lambda)\in U_{\Gamma}$.
 \item There is a smooth function
$z: S^1\times J_1\to J_2$ such that $\Gamma={\rm graph}\,z$.
\end{romanlist}
\end{Hypothesis}
As a consequence of \eqref{e:Gamma_rep} and Hypothesis~\ref{h:PtoE2}~(ii) we find
\begin{equation}\label{e:g_zeros}
 g(\varphi,\lambda_1,z(\varphi,\lambda_1))\equiv 0.
\end{equation}
As a transversality condition for $z$ we assume:
\begin{Hypothesis}\label{h:PtoE3}
 $z_\varphi(\varphi,0)=0 \quad\Rightarrow\quad  z_{\varphi\varphi}(\varphi,0)\not=0$.
\end{Hypothesis}
Fix some $\lambda_1^0$ close to zero, and let $\varphi^0$ be some value such that
$z_\varphi(\varphi^0,\lambda_1^0)=0$. Using this, we define
$\lambda_2^0:=z(\varphi^0,\lambda_1^0)$, and
$\lambda^0:=(\lambda_1^0,\lambda_2^0)$.  Now, considering the derivatives of
$g(\cdot,\lambda_1^0,z(\cdot,\lambda_1^0))$ at $\varphi=\varphi^0$ we find with
\eqref{e:g_zeros} and Hypothesis~\ref{h:PtoE3} that 
\begin{equation}\label{e:PtoE3_geom}
 g_\varphi(\varphi^0,\lambda^0)=0,\quad g_{\varphi\varphi}(\varphi^0,\lambda^0)\not=0.
\end{equation}
Note that ${\rm graph}\,g(\cdot,\lambda^0)$ describes the stable manifold
$W^s(E,\lambda^0)$ near $\varphi^0$. Therefore
\eqref{e:PtoE3_geom} means that $W^s(E,\lambda^0)$ and $W^u(P,\lambda^0)$ have a
quadratic tangency in $\varphi^0$.  We refer to
Figure~\ref{f:PtoE3_geom} for an illustration of the consequence of Hypothesis~\ref{h:PtoE3}.

Let the function $z(\cdot,\lambda_1)$ take its minimum in $\varphi_{min}(\lambda_1)$,
 and similarly let $z(\cdot,\lambda_1)$ be maximal in
$\varphi_{max}(\lambda_1)$.  This defines functions
$\lambda_{1,min}:\lambda_1\mapsto z(\varphi_{min}(\lambda_1),\lambda_1)$ and
$\lambda_{1,max}:\lambda_1\mapsto z(\varphi_{max}(\lambda_1),\lambda_1)$, both
mapping $J_1\to J_2$. The graphs of these functions define the $\lambda$-region
for which a heteroclinic cycle connecting $E$ and $P$ exists.  In our motivating
example this region is just the stripe between the two curves $t_b$
-- cf.  Figure~\ref{f:laser}.  Hence, the
maximal interval $[\underline{\lambda}_{\,2},\overline{\lambda}_2]$ defined in
Hypothesis~\ref{h:general} is given by
\begin{equation*}
 \underline{\lambda}_{\,2}:=z(\varphi_{min}(0),0), \quad \overline{\lambda}_2:=z(\varphi_{max}(0),0).
\end{equation*}
Moreover, for each $\lambda$ between the graphs of $\lambda_{1,min}$ and
$\lambda_{1,max}$, there are at least two heteroclinic PtoE
connection. These undergo saddle-node bifurcations on the graphs of
$\lambda_{1,min}$ and $\lambda_{1,max}$.  In particular,
moving along the $\lambda_2$-axis the heteroclinic PtoE connections undergo
saddle-node bifurcations in $\underline{\lambda}_{\,2}$ and
$\overline{\lambda}_2$.  If $z$ has exactly one minimum (and hence one
maximum), there are exactly two  heteroclinic PtoE connections
between the graphs of $\lambda_{1,min}$ and $\lambda_{1,max}$
-- cf.  Figure~\ref{f:graph_g}.

Next we consider the EtoP connection.
\begin{Hypothesis}\label{h:EtoP}
 There exist smooth functions $h^u: Q\to I$, $h^c: Q\to S^1$ such that 
\[
 \{(v^c,\delta,v^u)\in W^u(E,\lambda)\cap\Sigma^{in},\,\lambda\in Q\}
=\{(v^c,\delta,v^u)=(h^c(\lambda),\delta,h^u(\lambda)), \,\lambda\in Q\}.
\]
Moreover,
\begin{romanlist}
 \item $ h^u(0,\lambda_2)\equiv 0$,
 and\, $\forall\lambda_2\in J_2$ holds
$h^u_{\lambda_1}(0,\lambda_2)\not=0$,
\item $\exists q< 1:$ $\forall \varphi\in S^1\,\,|\frac{d}{d\varphi} h^c(0,z(\varphi,0))|\leq q$.
\end{romanlist}

\begin{figure}[ht]
\centering
\input{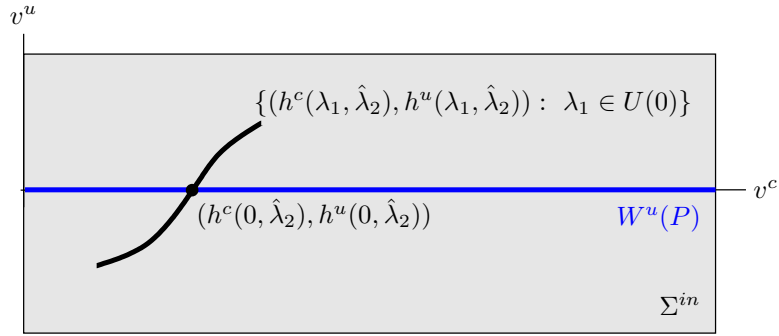}
 \caption{Visualisation of Hypothesis~\ref{h:EtoP}\,(i) with $\lambda_2=\hat\lambda_2$:
The curve $(h^c(\lambda_1,\hat\lambda_2),h^u(\lambda_1,\hat\lambda_2))=\mathop\bigcup\limits_{\lambda_1}W^u(E,\lambda_1,\hat\lambda_2)$ intersects $W^s(P)$ transversely.
}\label{f:graph_hu}
\end{figure}

\end{Hypothesis}
By definition
\begin{equation}\label{e:EtoP_char}
 W^u(E,\lambda)\cap  W^s(P,\lambda)\cap\Sigma^{in}\not=\emptyset\, \Longleftrightarrow\, h^u(\lambda)=0.
\end{equation}

So, Hypothesis~\ref{h:EtoP}\,(i) says that for all $\lambda$ on the
$\lambda_2$-axis, there is heteroclinic orbit connecting $E$
to $P$. In other words, the $\lambda_2$-axis is on a par with
the curve $c_b$ of our motivating example -- cf.
Figure~\ref{f:laser} or Figure~\ref{f:accum}, respectively.  Moreover, moving
through the $\lambda_2$-axis transversely effects that the EtoP connection
splits up with nonzero speed -- cf. also Lemma~\ref{l:accum}.
The consequences of Hypothesis~\ref{h:EtoP}\,(i) for the shape and mutual
position of the traces of $W^u(E)$ and $W^s(P)$ are depicted in
Figure~\ref{f:graph_hu}.  Finally, we note that by this
assumption $h^u_{\lambda_1}(0,\lambda_2)$ is bounded away from zero.

Recall that $z(\varphi,0)$ determines the $\lambda_2$ values for which a EtoP
cycle exists (clearly $\lambda_1=0$), where $\varphi$ is the
$v^c$-coordinate value of the intersection of the corresponding PtoE connection
with $\Sigma^{out}$. Whereas $h^c(0,z(\varphi,0))$ is the
$v^c$-coordinate value of the intersection of the corresponding EtoP connection
with $\Sigma^{in}$.  Hence, Hypothesis~\ref{h:EtoP}\,(ii)
yields that the proportion of the alteration rates of
these $v^c$-coordinates is bounded by $q<1$. In other
words, these coordinate values must not move against each
other too fast. Despite this geometric interpretation, this
hypothesis is more technical in nature. It will be used in the next section for
solving the bifurcation equations.

Now we can state our main result guaranteeing a snaking scenario.
\begin{Theorem}\label{t:snaking}
Assume Hypotheses~\ref{h:general} -- \ref{h:EtoP}. Then there is a constant $L_0>0$, and there are functions $\lambda_i:(L_0,\infty)\to \R$, $i=1,2$, such that
for each $L>L_0$ there is a one-homoclinic orbit to E for $\lambda\in Q$ that spends time $2L$ between $\Sigma^{in}$ and $\Sigma^{out}$ if and only if
$\lambda=(\lambda_1(L),\lambda_2(L)).$
\\[1ex]
Moreover there are an $\eta>0$, a $2\pi$-periodic function $\varphi_0^*(\cdot)$ and a positive bounded function $\hat{a}_u$ such that
\begin{align*}
 \lambda_1(L)&= {\textstyle{\frac{\hat{a}_u(L)}{h_{\lambda_1}^u(0,z(\varphi_0^*(L)+2L,0))}}}
e^{-2\alpha^u(0,z(\varphi_0^*(L)+2L,0))L} (1+ O(e^{-\eta L})),
\\[1ex]
 \lambda_2(L)&=z(\varphi^*_0(L)+2L,0)+O(e^{-\eta L}).
\end{align*}
\end{Theorem}

\begin{figure}[h]
\centering
\input{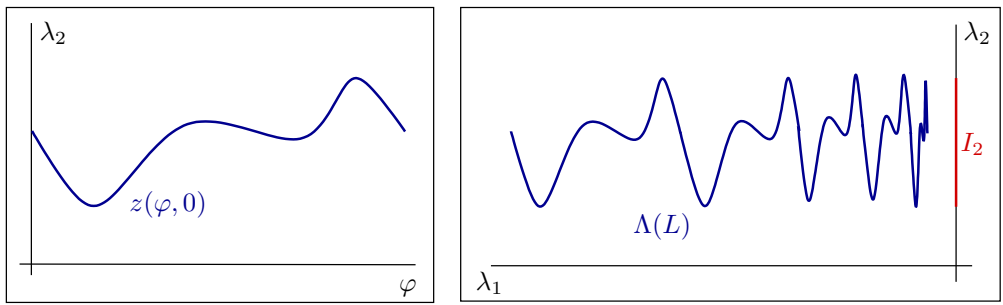}
\caption{The relation between ${\rm graph}\,z$ and the snaking curve $\Lambda(L)=(\lambda_1(L),\lambda_2(L))$. The shape of ${\rm graph}\, z(\cdot,0)$ depicted in the left panel is passed on to the snaking curve $\Lambda(L)$ in the right panel. The snaking curve accumulates at $\{0\}\times I_2$, the set of parameters for which the primary EtoP cycle exits, cf. Hypothesis~\ref{h:general} (iii).}\label{f:z_lambda_relat}
\end{figure}

It follows immediately that $\lambda_1(L)$ tends to zero as $L$
goes to infinity. Further it is obvious that $\lambda_2(L)$ is a perturbation of
$z(\varphi^*_0(L)+2L,0)$. This result resembles pretty much the statement about
the snaking parameter $\mu$ given in \cite[Theorem~2.2]{BKLSW:08}. But here, in
contrast to \cite{BKLSW:08}, the term $\varphi_0^*(L)$ is periodic and not
constant. If $\varphi^*_0(L)+2L$ is monotonically increasing, then
$\lambda_2(\cdot)$ essentially copies the behaviour of $z(\cdot,0)$, cf. Figure~\ref{f:z_lambda_relat} and Lemma~\ref{l:snaking_shape}.

In the
following lemma, we describe the shape of the snaking curve
$(\lambda_1(L),\lambda_2(L))$ somewhat closer. We consider $z(\cdot,0)$
as a periodic function $\R\to\R$. We denote the first and second derivative of $\lambda_2$
by $\lambda_2'$ and $\lambda_2''$, respectively.

\begin{Lemma}\label{l:snaking_shape}
  Assume Hypotheses~\ref{h:general} -- \ref{h:EtoP} with the more severe condition
$q<1/2$, cf. Hypothesis~\ref{h:EtoP}~(ii).
Then $\Phi:L\mapsto \varphi^*_0(L)+2L$ is a transformation, and for each $\hat\varphi$ with $z_\varphi(\hat\varphi,0)=0$ exists a unique $\hat L$ in a small neighbourhood of $\Phi^{-1}(\hat\varphi)$ such that $\lambda_2'(\hat L)=0$. Moreover $\lambda_2''(\hat L)\not=0$. These are the only zeros of $\lambda_2'$.
\end{Lemma}

\begin{figure}[ht]
\centering\input{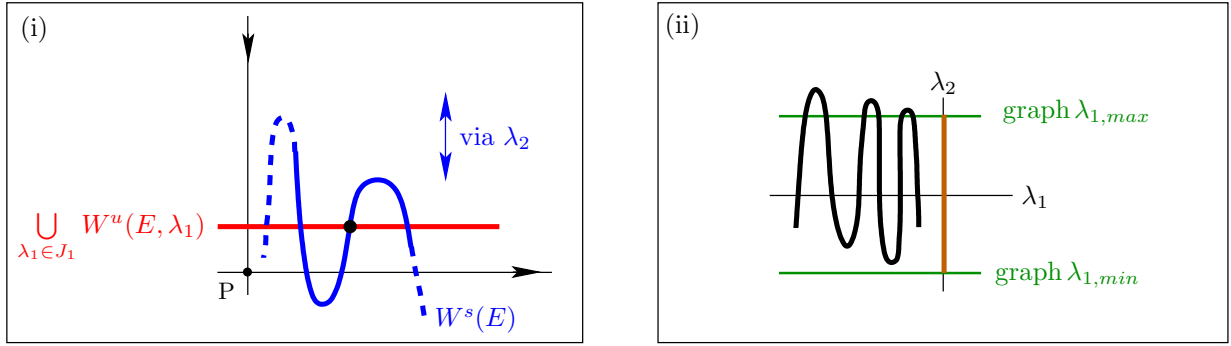}
 \caption{The creation of a snaking curve: Panel~(i) shows part of a (global) Poincar{\'e} section containing both $E$ and $P$. Panel~(ii) shows the continuation curve of one-homoclinic orbits to $E$.
}\label{f:geom_proof}
\end{figure}

The proofs of these statements are carried out in Section~\ref{s:snak_ana}.  Prior to that, however, we give
a geometric explanation with the help of the
Figure~\ref{f:geom_proof}.  Assume that the unstable
manifold of the equilibrium depends only on $\lambda_1$, and similarly that the
stable manifold of the equilibrium depends only on $\lambda_2$:
$W^u(E,\lambda)=W^u(E,\lambda_1)$,\quad $W^s(E,\lambda)=W^s(E,\lambda_2)$.  In
Figure~\ref{f:geom_proof}\,(i), we consider a fixed
Poincar{\'e} section of $P$. Fix some $\lambda_2$ -- and
therefore one particular position of $W^s(E)$ -- and assume
that an increasing $\lambda_2$ effects upward motion of $W^s(E)$. The bullet
defines a $\lambda_1$ for which $W^s(E)$ and $W^u(E)$ intersect and therefore a
homoclinic orbit to $E$ does exist. This homoclinic orbit can be continued by
moving $W^s(E)$ up and down. The corresponding continuation curve of the bullet
in the $\lambda$-space is displayed Figure~\ref{f:geom_proof}\,(ii).

\section{Snaking analysis}\label{s:snak_ana}
This section is devoted the proof of Theorem~\ref{t:snaking}.
A one-homoclinic orbit to $E$ can be conceived as built of three pieces: an orbit segment in $W^u(E)$ running from $E$ to $\Sigma^{in}$, a solution connecting $\Sigma^{in}$ and $\Sigma^{out}$ and an orbit segment in $W^s(E)$ running from $\Sigma^{out}$ to $E$.
Let $v$ be a solution starting in $\Sigma^{in}$ and arriving after time $2L$ in $\Sigma^{out}$. Then $v$ belongs to a one-homoclinic orbit to $E$ if the following coupling conditions are fulfilled
\begin{equation}\label{e:coupling_in_out}
 v(0,\lambda)\in W^u(E,\lambda)\cap \Sigma^{in},
\quad
v(2L,\lambda)\in W^s(E,\lambda)\cap \Sigma^{out}.
\end{equation}
As a consequence of the following Lemma~\ref{l:shilnikov} we get that for given $L$
there exists a unique solution $v$ starting in a certain
submanifold of $\Sigma^{in}$ and arriving after time $2L$ in $\Sigma^{out}$.
Afterwards, in the actual proof of Theorem~\ref{t:snaking} we use these
solutions to formulate coupling equations according to \eqref{e:coupling_in_out}.
The more general setting of Lemma~\ref{l:shilnikov} is used in Section~\ref{s:neg_Floqu_mul}.

\begin{Lemma}[Shilnikov problem near the periodic orbit]\label{l:shilnikov}
There is a positive constant $L_0$ such that for all $L>L_0$, all $(\varphi,\lambda)\in S^1\times Q$, and $\chi_s,\chi_u\in\{\pm 1\}$  
there exists a unique solution $v(t)$, also referred to as $v(t,\varphi,\lambda,\chi_s,\chi_u)$, of \eqref{e:fenichel} with
\[
 v^s(0)=\chi_s\delta,\quad v^c(0)=\varphi \quad \textrm{and} \quad v^u(2L)=\chi_u\delta.
\]
Moreover there is a positive constant $\eta < \min_{\lambda \in Q}\{|\alpha^s(\lambda)|, \alpha^u(\lambda) \}$ such that
\begin{equation}\label{e:v_est}
 \begin{array}{ll}
  v(0)&=(\varphi,\chi_s\delta, \chi_u a_u e^{-2\alpha^u(\lambda)L} \, (1+O(e^{-\eta L}))),
\\[1ex]
 v(2L)&=(\varphi+2L+O(e^{-\eta L}),\chi_s a_s e^{2\alpha^s(\lambda)L} \, (1+O(e^{-\eta L})),\chi_u\delta),
 \end{array}
\end{equation}
where $a_s$ and $a_u$ are positive functions depending on $(\varphi,\lambda,\chi_s,\chi_u)$. Moreover, 
$a_{s(u)}(\cdot,\cdot,\chi_s,\chi_u)$ are smooth.
For the derivatives of $v$ holds
\begin{equation}\label{e:dv_est}
 \begin{array}{ll}
  D_{\xi_1\ldots \xi_j} v(0)&=(D_{\xi_1\ldots \xi_j} \varphi,0, \chi_u D_{\xi_1\ldots \xi_j} (a_u  e^{-2\alpha^u(\lambda)L}) \, (1+O(e^{-\eta L}))),
\\[1ex]
 D_{\xi_1\ldots \xi_j} v(2L)&= (D_{\xi_1\ldots \xi_j}( \varphi+2L) + O(e^{-\eta L}), \chi_s  D_{\xi_1\ldots \xi_j}  (a_s e^{2\alpha^s(\lambda)L}) \, (1+O(e^{-\eta L})),0).
 \end{array}
\end{equation}
Here $\xi_i \in \{L, \lambda, \varphi\}$ for $i=1,\ldots, j$ and $j \in \{1,2,3\}$.
\end{Lemma}
Since for fixed $\chi_s, \chi_u$, the functions $a_s$ and $a_u$ are defined on the compact set $S^1\times Q$, they are bounded away from zero.

In \cite{Kru/San}, Krupa and Sandstede consider the Shilnikov problem in the context of slow/fast systems, where the slow manifold possesses a normally hyperbolic structure. Lemma \ref{l:shilnikov} is the counterpart to \cite[Theorem 4]{Kru/San}, and the proof proceeds in the same way as it is done there. 

\begin{prooof}[\,\,of Theorem~\ref{t:snaking}]

To describe the transition from $\Sigma^{in}$ to $\Sigma^{out}$ we use the function $v$ defined by Lemma~\ref{l:shilnikov}. By our choice of $\Sigma^{in}$ and $\Sigma^{out}$, cf. \eqref{e:def_Sigma}, this transition is determined by $v(\cdot,\varphi,\lambda):= v(\cdot,\varphi,\lambda,1,1)$.
Using the notation introduced in Section~\ref{s:setup_res}, equation \eqref{e:coupling_in_out} translates to
\begin{align}
 v^c(0,\varphi,\lambda)&=h^c(\lambda),\label{e:coup_in_hc}
\\
v^u(0,\varphi,\lambda)&=h^u(\lambda),\label{e:coup_in_hu}
\\
v^s(2L,\varphi,\lambda)&=g(v^c(2L,\varphi,\lambda),\lambda).\label{e:coup_out}
\end{align}

Recall that we are interested in those homoclinic orbits which are in a small
neighbourhood of a heteroclinic cycle. These cycles are
determined by $h^u=0$ and $g=0$, cf. \eqref{e:EtoP_char} and
\eqref{e:Gamma_rep}.  For that reason we solve
\eqref{e:coup_in_hc} -- \eqref{e:coup_out} near $h^u=0$ and $g=0$.

In accordance with \eqref{e:v_est} we find
that $v^c(2L,\varphi,\lambda)=\varphi+2L+O(e^{-\eta L})$.
Motivated by this equality, we introduce the following time transformation 
\begin{equation}\label{e:L_transf}
 2l=2L+O(e^{-\eta L}).
\end{equation}
Indeed, \eqref{e:L_transf} can be solved for 
\begin{equation}\label{e:L_transf_re}
 L=L_*(l,\lambda)=l+O(e^{-\eta l}).
\end{equation}
Using this, equation \eqref{e:coup_out} can be rewritten as $ v^s(2l+O(e^{-\eta l}),\lambda)=g(\varphi+2l,\lambda)$.
Altogether, using the new time $l$ and the estimates \eqref{e:v_est} the system \eqref{e:coup_in_hc} -- \eqref{e:coup_out} reads
\begin{align}
 \varphi&=h^c(\lambda_1,\lambda_2),\label{e:coup_in_hc_rew}
\\
a_u e^{-2\alpha^u(\lambda_1,\lambda_2)l}(1+O(e^{-\eta l}))&=h^u(\lambda_1,\lambda_2),\label{e:coup_in_hu_rew}
\\
a_s e^{2\alpha^s(\lambda_1,\lambda_2)l}(1+O(e^{-\eta l}))&=g(\varphi+2l,\lambda_1,\lambda_2).\label{e:coup_out_rewr}
\end{align}
First, we consider \eqref{e:coup_out_rewr}, which describes the coupling in $\Sigma^{out}$.
Recall that we want to solve \eqref{e:coup_out}, and therefore also \eqref{e:coup_out_rewr}, near $g=0$. Furthermore, recall from \eqref{e:g_zeros} that 
$g(\varphi+2l,\lambda_1,z(\varphi+2l,\lambda_1))\equiv 0$.
Now
write 
\begin{equation*}
\lambda_2=z(\varphi+2l,\lambda_1)+\mu 
\end{equation*}

and expand $g(\varphi+2l,\lambda_1,z(\varphi+2l,\lambda_1)+\mu)$ w.r.t. $\mu$. Inserting in \eqref{e:coup_out_rewr} gives
 \begin{equation*}
 a_s e^{2\alpha^s(\lambda_1,z(\varphi+2l,\lambda_1)+\mu)l}(1+O(e^{-\eta l}))
=g_{\lambda_2}(\varphi+2l,\lambda_1,z(\varphi+2l,\lambda_1))\mu+O(\mu^2).
\end{equation*}
Using this and Hypothesis~\ref{h:PtoE2}, the coupling equation \eqref{e:coup_out} eventually reads:
\begin{equation}\label{e:coup_out_rewr2}
{\textstyle{
 \frac{a_s(\varphi,\lambda_1,z(\varphi+2l,\lambda_1)+\mu)}{g_{\lambda_2}(\varphi+2l,\lambda_1,z(\varphi+2l,\lambda_1))}}}
e^{2\alpha^s(\lambda_1,z(\varphi+2l,\lambda_1)+\mu)l}(1+O(e^{-\eta l}))
=\mu+O(\mu^2).
\end{equation}
For $|\mu|\ll 1$, sufficiently large $l$ and all $\varphi$ this equation can be solved for $\mu=\mu^*(l,\varphi,\lambda_1)$ by means of the implicit function theorem \cite[Chapter~2.2]{ChowHale82}.
The solving function $\mu^*$ is differentiable.
Further we see from \eqref{e:coup_out_rewr2} that
\begin{equation}\label{e:mu*_rep}
 \mu^*(l,\varphi,\lambda_1)=
{\textstyle{
 \frac{a_s(\varphi,\lambda_1,z(\varphi+2l,\lambda_1))}{g_{\lambda_2}(\varphi+2l,\lambda_1,z(\varphi+2l,\lambda_1))}}}
e^{2\alpha^s(\lambda_1,z(\varphi+2l,\lambda_1))l}(1+O(e^{-\eta
l})).
\end{equation}
Altogether we find that the coupling equation \eqref{e:coup_out} can be solved for $\lambda_2=\hat\lambda_2(l,\varphi,\lambda_1)$ with
\begin{equation*}\label{e:lambda2*_rep}
 \hat\lambda_2(l,\varphi,\lambda_1)=z(\varphi+2l,\lambda_1))+\mu^*(l,\varphi,\lambda_1),
\end{equation*}
where the leading order term of $\mu^*$ is given by \eqref{e:mu*_rep}.

Now we turn towards the coupling in $\Sigma^{in}$ which is determined by \eqref{e:coup_in_hc_rew} and \eqref{e:coup_in_hu_rew}.
Using the representation of $\hat\lambda_2$, these equations read
\begin{align}
 \varphi&= h^c(\lambda_1,z(\varphi+2l,\lambda_1)+\mu^*(l,\varphi,\lambda_1)),\label{e:coup_in_hc_rewr}
\\
a_u e^{-2\alpha^u(\lambda_1,z(\varphi+2l,\lambda_1)+\mu^*(l,\varphi,\lambda_1))l}(1+O(e^{-\eta l}))
&= h^u(\lambda_1,z(\varphi+2l,\lambda_1)+\mu^*(l,\varphi,\lambda_1)).\label{e:coup_in_hu_rewr}
\end{align}
We solve \eqref{e:coup_in_hc_rewr}, \eqref{e:coup_in_hu_rewr}
for $(\varphi,\lambda_1)$ depending on $l$. Note that $\varphi$ is the $v^c$
coordinate where the prospective homoclinic orbits hits $\Sigma^{in}$.
Hence, $\varphi$ may vary within a ``large'' range. To handle
this difficulty analytically, we first consider the
``unperturbed equation'' $\varphi=h^c(0,z(\varphi+2l,0))$.  For that we consider
$z(\cdot,0)$ as a $2\pi$-periodic function $\R\to\R$, cf.
Hypotheses~\ref{h:general}~and~\ref{h:PtoE2}. Hence, 
$h^c(0,z(\cdot+2l,0))$ is a $2\pi$-periodic function as well.  Because of
Hypothesis~\ref{h:EtoP}\,(ii), we can apply again the implicit
function theorem to find a unique solution $\varphi^*_0(l)$ on $\R$ such that
\begin{equation}\label{e:coup_in_hc_rewr_unp}
 \varphi=h^c(0,z(\varphi+2l,0))\,\Longleftrightarrow\,\varphi=\varphi^*_0(l).
\end{equation}
Note that $\varphi^*_0(\cdot)$ is again $2\pi$-periodic.
Now, write 
\begin{equation*}
 \varphi=\varphi^*_0(l)+\psi,
\end{equation*}

and we define $H^c(l,\psi,\lambda_1)$, $H^u(l,\psi,\lambda_1)$ by 
\begin{equation*}
 H^{c/u}(l,\psi,\lambda_1):=h^{c/u}(\lambda_1,z(\varphi^*_0(l)+\psi+2l,\lambda_1)).
\end{equation*}
Using these terms the right-hand sides of \eqref{e:coup_in_hc_rewr} and \eqref{e:coup_in_hu_rewr} read
\begin{equation*}
h^{c/u}(\lambda_1,z(\varphi^*_0(l)+\psi+2l,\lambda_1)+\mu^*(l,\varphi^*_0(l)+\psi,\lambda_1))
=H^{c/u}(l,\psi,\lambda_1)+r^{c/u}(\mu^*(l,\varphi^*_0(l)+\psi,\lambda_1)),
\end{equation*}
where $r^{c/u}(\mu^*)=O(\mu^*)$.
Further, since $\varphi^*_0$ is the unique solution of
\eqref{e:coup_in_hc_rewr_unp}, we find
\begin{equation*}\label{e:Hc_expand}
 H^c(l,0,0)=h^c(0,z(\varphi^*_0(l)+2l,0))=\varphi^*_0(l).
\end{equation*}
In accordance with Hypothesis~\ref{h:EtoP}\,(i), we find furthermore
\begin{equation*}
H^u(l,0,0)=0,\quad H^u_\psi(l,0,0)=0.
\end{equation*}
Hence, \eqref{e:coup_in_hc_rewr} and \eqref{e:coup_in_hu_rewr} are equivalent to 
\begin{align}
\psi&=H^c_\psi(l,0,0)\psi + H^c_{\lambda_1}(l,0,0)\lambda_1 + O(|(\psi,\lambda_1)|^2)+r^c(\mu^*)\label{e:coup_in_hc_rewr3},
\\
a_u e^{-2\hat\alpha^u(l,\psi,\lambda_1)l}(1+O(e^{-\eta l}))
&=H^u_{\lambda_1}(l,0,0)\lambda_1 + O(|(\psi,\lambda_1)|^2)+r^u(\mu^*),\label{e:coup_in_hu_rewr3}
\end{align}
with $\hat\alpha^u(l,\psi,\lambda_1):= \alpha^u(\lambda_1,z(\varphi^*_0(l)+\psi+2l,\lambda_1)+\mu^*(l,\varphi^*_0(l)+\psi,\lambda_1))$.

Our goal is now to solve the system \eqref{e:coup_in_hc_rewr3},
\eqref{e:coup_in_hu_rewr3} for $(\psi,\lambda_1)$ depending on $l$.
To this end we invoke again the implicit
function theorem.  The main observation is that due to Hypothesis~\ref{h:EtoP} 
\begin{equation*}
 |H^c_\psi(l,0,0)|\leq q< 1,\quad  H^u_{\lambda_1}(l,0,0)\not=0.
\end{equation*}
Note that due to \eqref{e:dv_est}, the corresponding partial
derivatives of $r^c(\mu^*(l,\varphi^*_0(l)+\psi,\lambda_1))$ and
$r^u(\mu^*(l,\varphi^*_0(l)+\psi,\lambda_1))$ tend to zero as $l\to\infty$.
Therefore there exist unique functions $\psi^*(l)$,
$\lambda_1^*(l)$ satisfying the system \eqref{e:coup_in_hc_rewr3},
\eqref{e:coup_in_hu_rewr3}.  Accordingly,
\eqref{e:coup_in_hc_rewr}, \eqref{e:coup_in_hu_rewr} are satisfied by 
\begin{equation*}
 \phi^*(l):=\varphi^*_0(l)+\psi^*(l)
\quad{\rm and}\quad
\lambda_1^*(l).
\end{equation*}

Inspecting \eqref{e:coup_in_hu_rewr} and
\eqref{e:coup_in_hc_rewr3}, we find with 
$\hat a_u(l):=a_u(\varphi_0^*(l),0,z(\varphi_0^*(l),0))$ that
\begin{equation*}
 \lambda_1^*(l)= {\textstyle{\frac{\hat a_u(l)}{H_{\lambda_1}^u(l,0,0)}}}
e^{-2\hat\alpha^u(l,0,0)l} (1+ O(e^{-\eta l}))
\quad{\rm and}\quad
\psi^*(l)=O(e^{-2\hat\alpha^u(l,0,0) l}).
\end{equation*}

Altogether, for \eqref{e:coup_in_hc_rew} -- \eqref{e:coup_out_rewr} we find the unique solution 
$(\varphi,\lambda_1,\lambda_2)(l)=(\phi^*(l),\lambda_1^*(l),\lambda_2^*(l))$, where
\begin{align*}
 \lambda_2^*(l):=\hat\lambda_2(l,\phi^*(l),\lambda_1^*(l))
&=z(\phi^*(l)+2l,\lambda_1^*(l))+\mu^*(l,\phi^*(l),\lambda_1^*(l))
\\
&=z(\varphi^*_0(l)+2l,0)+O(e^{-\eta l}).
\end{align*}

Note that $v$ spends time  $2l+O(e^{-\eta l})$ between $\Sigma^{in}$ and $\Sigma^{out}$, cf. \eqref{e:L_transf}.
So, in view of the statement in Theorem~\ref{t:snaking}, we define
\begin{equation*}
 \phi(L):=\phi^*(l(L)),\quad \lambda_1(L):=\lambda_1^*(l(L)),\quad \lambda_2(L):=\lambda_2^*(l(L)).
\end{equation*}
Then $(\varphi,\lambda_1,\lambda_2)(L)=(\phi(L),\lambda_1(L),\lambda_2(L))$ solves  \eqref{e:coup_in_hc} -- \eqref{e:coup_out}.
The above considerations yield
\begin{align*}\label{e:phi(L)_rep}
 \phi(L)&=\varphi_0^*(L)+O(e^{-\eta L}),
\\[1ex]
 \lambda_1(L)&= {\textstyle{\frac{\hat a_u(L)}{H_{\lambda_1}^u(L,0,0)}}}
e^{-2\hat\alpha^u(L,0,0)L} (1+ O(e^{-\eta L})),
\\[1ex]
 \lambda_2(L)&=z(\varphi^*_0(L)+2L,0)+O(e^{-\eta L}).
\end{align*}
This finally completes the proof of Theorem~\ref{t:snaking}.
\end{prooof}

\begin{prooof}[\,\,of Lemma~\ref{l:snaking_shape}]

In what follows, we sketch the proof of Lemma~\ref{l:snaking_shape}. 
We note that, due to \eqref{e:dv_est}, the $O$-term in the representation of $\lambda_2(L)$ is differentiable and its derivative can be estimated by a $O$-term of the same order. The same holds true for higher derivatives.
Therefore we find
\begin{equation*}\label{e:lambda2_der}
 \lambda_2'(L)=z_\varphi(\varphi_0^*(L)+2L,0)({\varphi_0^*}'(L)+2)
+O(e^{-\eta L}).
\end{equation*}
Recall the determining equation $\varphi_0^*(L)=h^c(0,z(\varphi_0^*(L)+2l,0))$
for $\varphi_0^*(L)$, cf. \eqref{e:coup_in_hc_rewr_unp}.
From this equation, we get an estimate of the derivative of
$\varphi_0^*(L)$,
whereby we finally confirm that for sufficiently large $L$
\begin{equation*}
{\varphi_0^*}'(L)+2\not=0.
\end{equation*}
So, necessarily the zeros of $\lambda_2'(L)$ are close to the zeros of
$z_\varphi(\varphi_0^*(L)+2L,0))$.  Let $z_\varphi(\varphi_0^*(L_0)+2L_0,0))=0$.
Using the contraction principle we find a neighbourhood $U(L_0)$ of $L_0$ in
which $\lambda_2'(L)=0$ has a unique solution $\hat L$.  Straightforward
computations show $\lambda_2''(\hat L)\not=0$.  Further, the size of the
neighbourhood $U(L_0)$ can be chosen independently of $L_0$. Outside the union
of these neighbourhoods, $z_\varphi(\cdot,0)$ is bounded away
from zero. This finally shows that outside the union of these neighbourhoods
$\lambda_2'(L)$ has no zeros for sufficiently large $L$.
\end{prooof}

\section{Negative Floquet multipliers}\label{s:neg_Floqu_mul}
Now we discuss the scenario, where the nontrivial Floquet multipliers of the periodic orbit $P$ are negative -- in other words, we replace Hypothesis \ref{h:general} by: 
\begin{Hypothesis}\label{h:neg_Floquet_mul}
 The nontrivial Floquet multipliers of $P$ are negative.
\end{Hypothesis}
Recall, since our setting is in $\R^3$, the two nontrivial Floquet multipliers of $P$ must have the same sign. 
Negative multipliers cause that the vector bundle consisting of the eigenvectors of the monodromy matrices along $P$ is a M\"obius strip and thus not orientable. Hence, we cannot introduce Fenichel coordinates near $P$.
We overcome this difficulty by introducing local coordinates, which are not $2 \pi$ periodic, but of period $4 \pi$. To this end we transform at first \eqref{e:system} into normal form, cf. \cite[Theorem 3.11]{SSTC98}, which gives:
\begin{align*}
\dot{\theta}&=1, \cr
\dot{y}& = B(\theta)y + F(\theta,y,\lambda),
\end{align*}
where $y=(y_1,y_2) \in \R^2$ and $\theta \in S^1$; $B$ and $F$ are smooth.
Furthermore we straighten the stable and unstable fibres of $P$, as it is done in \cite{Jon95}. This yields that the function $F$ satisfies
\begin{equation*}
F(\theta,0,y_2,\lambda) \equiv 0 \quad \textrm{and} \quad F(\theta,y_1,0,\lambda)\equiv 0 \quad \textrm{and} \quad D_yF(\theta,0,0,\lambda)\equiv 0, \quad \forall \, \lambda \in Q, \, \forall \theta \in S^1.
\end{equation*}
After that we apply Floquet theory to the linear system $\dot{y} = B(\theta)y$, see \cite[Theorem 3.12]{SSTC98} for more details. This transforms the above normal form into:
\begin{align}
\dot{v^c}&=1,\nonumber\\
\dot{v^s}& =  {\alpha^s(\lambda)} v^s + \tilde{F}^s(v^c,v^s,v^u,\lambda),\label{e:floquet} \\
\dot{v^u}& =  \alpha^u(\lambda) v^u + \tilde{F}^u(v^c,v^s,v^u,\lambda),\nonumber
\end{align}
where $v^c \in {\mathbb S}^1:=\R/_{\sim_{4\pi}}$, and $x\sim_{4\pi} y \Leftrightarrow x=y\!\mod\!4\pi$.
Moreover, similar arguments as leading to Lemma~\ref{l:fenichel} yield $\tilde{F}^s(v^c,0,v^u,\lambda)=\tilde{F}^u(v^c,v^s,0,\lambda)= 0$, as well as, $D_{v^s}\tilde{F}^s(v^c,0,0,\lambda)=D_{v^u}\tilde{F}^u(v^c,0,0,\lambda)= 0$.

Note that by this construction the two points $(v^c,v^s,v^u)$ and $(v^c+2\pi,-v^s,-v^u)$ represent the same point in $(\theta,y)$-coordinates.
In other words, 
two points are identified via the map:
\begin{align*}\label{e:ident}
i  : \,  {\mathbb S}^1 \times I \times I & \rightarrow  {\mathbb S}^1 \times I \times I \cr
(v^c,v^s,v^u) & \mapsto (v^c+2\pi,-v^s,-v^u).
\end{align*}
Next, in accordance with the procedure in Section~\ref{s:setup_res}, we introduce a cross-section $\Sigma^{in}$ of $W^s(P)$ intersecting orthogonally the stable fibres of $P$ in a distance $\delta$ of $P$. Similarly we define $\Sigma^{out}$.
In $(v^c,v^s,v^u)$-coordinates these sections read:
\begin{equation}\label{e:def_Sigma-}
 {\Sigma}^{in}_{+}:={\mathbb S}^1 \times \{v^s=\delta\} \times I, \qquad {\Sigma}^{out}_{+}:={\mathbb S}^1 \times I \times \{v^u=\delta\}.
\end{equation}
The subscript ``+'' refers to the positive value $\delta$ for the fixed $v^s$- and $v^u$-coordinate, respectively.
The sections defined in \eqref{e:def_Sigma-} are identified via the map $i$ with
\begin{equation*}\label{e:def_Sigma+}
 {\Sigma}^{in}_{-}:={\mathbb S}^1 \times \{v^s=-\delta\} \times I, \qquad {\Sigma}^{out}_{-}:={\mathbb S}^1 \times I \times \{v^u=-\delta\}.
\end{equation*}
Further we introduce a set $\Gamma$ similarly to \eqref{e:Gamma_def} -- formally replacing $S^1$ by ${\mathbb S}^1$:
\begin{equation}\label{e:Gamma_def_neg}
 {\Gamma}:=\{(\varphi,\lambda)\in {\mathbb S}^1\times Q: W^s(E,\lambda)\cap W^{uu}(P(\varphi,\lambda),\lambda)\cap {\Sigma}^{out}_{+}\not=\emptyset\}.
\end{equation}
Note that each strong unstable fibre $W^{uu}(P(\varphi,\lambda),\lambda)$ of $P$ intersects $\Sigma^{out}$ twice. In the terminology of \eqref{e:Gamma_def_neg} those two points are represented by $W^{uu}(P(\varphi,\lambda),\lambda)\cap {\Sigma}^{out}_{+}$ and $W^{uu}(P(\varphi+2\pi,\lambda),\lambda)\cap {\Sigma}^{out}_{+}$.

Let $U_\Gamma$ be an open neighbourhood of $\Gamma$ in ${\mathbb S}^1\times Q$. Regarding the PtoE connecting orbit we assume:
\begin{Hypothesis}\label{h:PtoE1_neg}
 There is a smooth function $g:U_{\Gamma}\to I$ and an $\epsilon >0$ such that
\begin{equation*}
  \{(\varphi,v^s,\delta)\in W^s(E,\lambda)\cap\Sigma^{out}_{+}: |v^s|<\epsilon, (\varphi,\lambda)\in U_{\Gamma} \}
=\{(\varphi,g(\varphi,\lambda),\delta): (\varphi,\lambda)\in U_{\Gamma} \}.
\end{equation*}
\end{Hypothesis}

\begin{Hypothesis}\label{h:PtoE2_neg}
 \begin{romanlist}
 \item There is a  constant $b>0$ such that $|g_{\lambda_2}(\varphi,\lambda)|\geq b$, for all $(\varphi,\lambda)\in U_{\Gamma}$.
 \item There is a smooth function
$z: {\mathbb S}^1\times J_1\to J_2$ such that $\Gamma={\rm graph}\,z$.
\end{romanlist}
\end{Hypothesis}
Consequently
\begin{equation*}
 g(\varphi,\lambda_1,z(\varphi,\lambda_1))\equiv 0.
\end{equation*}

\begin{Hypothesis}\label{h:PtoE3_neg}
 $z_{\varphi}(\varphi,0)=0 \quad\Rightarrow\quad  z_{\varphi\varphi}(\varphi,0)\not=0$.
\end{Hypothesis}

Finally, regarding the EtoP connecting orbit we assume:
\begin{Hypothesis}\label{h:EtoP_neg}
 There exist smooth functions $h^u: Q\to I$, $h^c: Q\to {\mathbb S}^1$ such that
\[
 \{(v^c,\delta,v^u)\in W^u(E,\lambda)\cap\Sigma_{+}^{in},\,\lambda\in Q\}
=\{(v^c,\delta,v^u)=(h^c(\lambda),\delta,h^u(\lambda)), \,\lambda\in Q\}.
\]
Moreover,
\begin{romanlist}
 \item $ h^u(0,\lambda_2)\equiv 0$,
 and\, $\forall\lambda_2\in J_2$  holds
$h^u_{\lambda_1}(0,\lambda_2)>0$,
\item $\exists q< 1:$ $\forall \varphi\in {\mathbb S}^1\,\,|\frac{d}{d\varphi} h^c(0,z(\varphi,0))|\leq q$.
\end{romanlist}
\end{Hypothesis}
Indeed, in Hypothesis~\ref{h:EtoP_neg}(i) it already suffices to assume $h^u_{\lambda_1}(0,\lambda_2)\not=0$. The specification stated in the hypothesis determines the sign of the functions $\lambda_1^+$ and $\lambda_1^-$ in the way as stated in the theorem below.

Now, the analogue of Theorem~\ref{t:snaking} reads:
\begin{Theorem}\label{t:snaking_floquet}
Assume Hypothesis~\ref{h:general}, and Hypotheses~\ref{h:neg_Floquet_mul} -- \ref{h:EtoP_neg}.
Then there is a constant $L_0>0$, and there are functions $\lambda_i^{+},\lambda_i^{-} :(L_0,\infty)\to \R$, $i=1,2$, such that
for each $L>L_0$ there is a one-homoclinic orbit to E for $\lambda\in Q$ that spends time $2L$ between $\Sigma^{in}$ and $\Sigma^{out}$ if and only if
$\lambda=(\lambda_1^\pm(L),\lambda_2^\pm(L)).$
\\[1ex]
Moreover there are an $\eta>0$, two $4\pi$-periodic functions $\varphi_0^+(\cdot)$, $\varphi_0^-(\cdot)$ and positive bounded functions $\hat{a}_u^+$, $\hat{a}_u^-$  such that
\begin{align*}
 \lambda_1^\pm(L)&= {\textstyle{\frac{\pm\hat{a}_u^{\pm}(L)}{ h_{\lambda_1}^u(0,z(\varphi_0^\pm(L)+2L,0))}}}
e^{-2{\alpha^u(0,z(\varphi_0^\pm(L)+2L,0))}L} (1+ O(e^{-\eta L})),
\\[1ex]
 \lambda_2^\pm(L)&=z(\varphi^\pm_0(L)+2L,0)+O(e^{-\eta L}). \cr
\end{align*}
\end{Theorem}
A visualisation of the statement of this theorem is given in Figure~\ref{f:z_lambda_relat_neg}.
\begin{figure}[ht]
\centering
\input{figures/neg-mult.pspdftex}
\caption{As in the case of positive multipliers the shape of ${\rm graph}\, z(\cdot,0)$ is passed on to the snaking curves $\Lambda^\pm(L)$, cf. Figure~\ref{f:z_lambda_relat}.
The snaking curves accumulate at $\{0\}\times I_2$ from different sides. }\label{f:z_lambda_relat_neg}
\end{figure}
\begin{prooof}
We pursue the same strategy as in the proof of Theorem~\ref{t:snaking}: We construct one-homoclinic orbits to $E$ by coupling in $\Sigma^{in}$ and $\Sigma^{out}$ the unstable and stable manifolds, respectively, with solutions according to Lemma~\ref{l:shilnikov}. However, here there are two possibilities for the transition from $\Sigma^{in}$ to $\Sigma^{out}$. There are solutions of \eqref{e:floquet} starting in $\Sigma_+^{in}\cap\{v^u>0\}$ and end up in $\Sigma_+^{out}\cap\{v^s>0\}$, and there are solutions starting in $\Sigma_-^{in}\cap\{v^u>0\}$ and end up in $\Sigma_+^{out}\cap\{v^s<0\}$, cf. Figure~\ref{f:sigmas_neg}. In the language of Lemma~\ref{l:shilnikov} this distinction is determined by the signs of $\chi_s$ and $\chi_u$.
\begin{figure}[ht]
\centering
\input{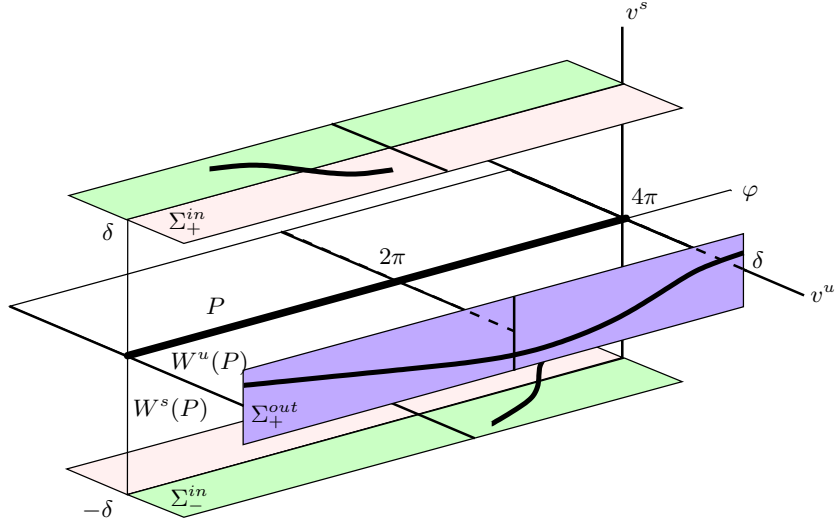}
\caption{The cross-sections $\Sigma^{in}$ and $\Sigma^{out}$. In $\Sigma^{in}_+$ and $\Sigma^{in}_-$ there is drawn a curve $(h^c,h^u)(\cdot,\lambda_2)$ for fixed $\lambda_2$, cf. Hypothesis~\ref{h:EtoP_neg}. In $\Sigma^{out}_+$ there is drawn the graph of $g(\cdot,\lambda)$ for fixed $\lambda$, cf. Hypothesis\ref{h:PtoE1_neg}.}\label{f:sigmas_neg}
\end{figure}

{\bf The transition $\Sigma_+^{in}\cap\{v^u>0\}$ to $\Sigma_+^{out}\cap\{v^s>0\}$}

Here we employ solutions of Lemma~\ref{l:shilnikov} with $(\chi_s,\chi_u)=(1,1)$ what we suppress from the notation.
In this case the argumentation runs completely parallel to the proof of Theorem~\ref{t:snaking}. We confine to sketch the procedure:
The coupling equations analogue to \eqref{e:coup_in_hc}-\eqref{e:coup_out} are almost the same:
 \begin{align*}
v^c(0,\varphi,\lambda)& = h^c(\lambda),\nonumber\\
v^u(0,\varphi,\lambda) & = h^u(\lambda),\label{e:coup-q1}\\
v^s(2L,\varphi,\lambda) & =g(v^c(2L,\varphi,\lambda),\lambda),\nonumber
\end{align*}
with the only difference that here $\varphi, v^c\in{\mathbb S}^1$. 
From that we gain the analogue to \eqref{e:coup_in_hc_rew}-\eqref{e:coup_out_rewr}
\begin{eqnarray}
  \varphi&\!\!\!\!=&\!\!\!\!h^c(\lambda),\nonumber
\\
   a_u e^{-2{\alpha^u(\lambda)}l}(1+O(e^{-\eta l}))&\!\!\!\!=&\!\!\!\! h^u(\lambda),\label{e:hu+}
\\
  a_s e^{2{\alpha^s(\lambda)}l}(1+O(e^{-\eta l}))&\!\!\!\!=&\!\!\!\! g(\varphi+2l,\lambda).\nonumber
\end{eqnarray}
Proceeding in the same way as in the proof of Theorem~\ref{t:snaking} we get the solutions $\lambda_1^+(L)$ and $\lambda_2^+(L)$ as stated in the theorem. The function $\varphi_0^+$ solves the analogue to \eqref{e:coup_in_hc_rewr_unp}.

{\bf The transition $\Sigma_-^{in}\cap\{v^u>0\}$ to $\Sigma_+^{out}\cap\{v^s<0\}$}

First we note that the intersections of $W^u(E,\lambda)$ with $\Sigma^{in}$ written in the form
$(h^c(\lambda),\delta,h^u(\lambda))\in \Sigma_+^{in}$, cf. Hypothesis~\ref{h:EtoP_neg},
are identified with 
$(h^c(\lambda)+2\pi,-\delta,-h^u(\lambda))\in \Sigma_-^{in}$.
This allows to employ solutions of Lemma~\ref{l:shilnikov} with $(\chi_s,\chi_u)=(-1,1)$ for our analysis. Again we suppress the $\chi_s,\chi_u$-dependence from the notation. 

Thus the coupling equations read:
\begin{align*}
v^c(0,\varphi,\lambda) & =h^c(\lambda) +2\pi,\cr
v^u(0,\varphi,\lambda)&= - h^u(\lambda),
\\
v^s(2L,\varphi,\lambda)&= g(v^c(2L,\varphi,\lambda),\lambda).\notag
\end{align*}
With the results of Lemma~\ref{l:shilnikov} this can be rewritten as
\begin{eqnarray}
  \varphi&\!\!\!\!=&\!\!\!\! h^c(\lambda)+2\pi,\nonumber
\\
  a_u e^{-2{\alpha^u(\lambda)}l}(1+O(e^{-\eta l}))&\!\!\!\!=&\!\!\!\! -h^u(\lambda),\label{e:hu-}
\\
- a_s e^{2{\alpha^s(\lambda)}l}(1+O(e^{-\eta l}))&\!\!\!\!=&\!\!\!\! g(\varphi+2l,\lambda).\nonumber
\end{eqnarray}
Now we can proceed again as in the proof of Theorem~\ref{t:snaking}, with the minor difference that the function $\varphi_0^-$ results from the fixed point equation
\begin{equation*}
 \varphi=h^c(0, z(\varphi+2l,0))+2\pi.
\end{equation*}
This finally leads to the solutions $\lambda_1^-(L)$ and $\lambda_2^-(L)$.

We want to note that $\lambda_1^+$ and $\lambda_1^-$ have different signs. This follows immediately from \eqref{e:hu+} and \eqref{e:hu-} in combination with Hypothesis~\ref{h:EtoP_neg}(i).
\end{prooof}

\section{A nonsnaking scenario}\label{s:nonsnak_ana}
Now we consider a further continuation scenario of homoclinic orbits.  We retain
the general Hypothesis~\ref{h:general}, and assume positive Floquet multipliers. Further we adopt the definition of
$\Gamma$, cf. \eqref{e:Gamma_def}, and assume Hypothesis~\ref{h:PtoE1}.  In
Section \ref{s:snak_ana}, we have seen that
the shape of the continuation curve $(\lambda_1(L),\lambda_2(L))$
is basically determined by the form of the set $\Gamma$, which was assumed to be
the graph of a function $z \, : \, S^1 \times J_1 \rightarrow J_2$. Now we
investigate the consequences of altering the corresponding
Hypothesis~\ref{h:PtoE2}. More precisely, we  suppose that $\Gamma$ is no longer
graph of a function, but for fixed $\lambda_1$ a closed curve in each case. 

In contrast to the situation in Section~\ref{s:snak_ana}, we end up with a sequence of closed
continuation curves, so-called isolas.  Consequently, the
one-homoclinic orbits close to the primary EtoP cycle do not lie on one (global)
continuation curve.  The addressed curves tend to $\{0\} \times I_2$ in the
sense of the Hausdorff metric . 

A similar scenario was already discussed in \cite{BKLSW:08} for the Hamiltonian case under slightly different assumptions.

We assume the following:

\begin{Hypothesis}\label{h:PtoE2-non-snake}
Let $I_1 \subsetneq J_1$ be a closed  interval containing zero and let $I_\varphi \subsetneq S^1$. There exist smooth functions $\hat{\varphi} \, : \, S^1 \times I_1 \rightarrow I_\varphi$ and  $\hat{\lambda}_2  \, : \, S^1 \times I_1 \rightarrow J_2 $ 
such that
\begin{equation*}
 \Gamma  = \{ (\hat{\varphi}(r, \lambda_1), \, \lambda_1,  \, \hat{\lambda}_2(r, \lambda_1)) \, : \, r \in S^1, \, \lambda_1 \in I_1  \},
\end{equation*}
where
\begin{equation*} 
(D_r\hat{\varphi}(r, \lambda_1), D_r\hat{\lambda}_2(r, \lambda_1)) \neq 0, \qquad \forall r \in S^1, \,\forall \, \lambda_1 \in I_1.
\end{equation*}
\end{Hypothesis}
Hence, the set $\Gamma\cap\{\lambda_1={\rm const.}\}$ is a closed, regular curve in $I_\varphi \times J_2$ parametrised by some parameter $r \in S^1$. 

\begin{figure}[ht]
\centering
\input{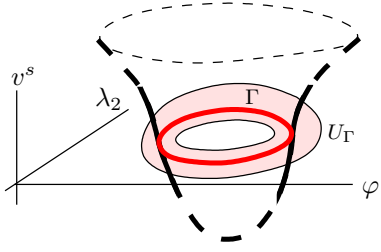}
 \caption{Visualisation of Hypothesis~\ref{h:PtoE2-non-snake}. Depicted is the closed curve 
$\{ (\hat{\varphi}(r, 0), \, \lambda_1,  \, \hat{\lambda}_2(r,0)) \, : \, r \in S^1\} = \{g(\cdot,0,\cdot)=0\}$.}
\end{figure}

Again we denote by  $U_{\Gamma}$ an open neighbourhood of $\Gamma$ in $S^1
\times Q$. Further, \eqref{e:Gamma_rep} together with
Hypothesis \ref{h:PtoE2-non-snake} yield as counterpart to \eqref{e:g_zeros}
\begin{equation}\label{e:Gamma_rep-non-snake}
g(\hat{\varphi}(r, \lambda_1), \lambda_1,\hat{\lambda}_2(r, \lambda_1) ) \equiv 0, \qquad \forall \, (r,\lambda_1)\in S^1 \times I_1.
\end{equation}

\begin{Hypothesis}\label{h:unique}
\begin{equation*}
D_{\varphi}	g(\hat{\varphi},0,\hat{\lambda}_2) D_r\hat{\lambda}_2(r, 0) - D_{\lambda_2}	g(\hat{\varphi},0,\hat{\lambda}_2) D_r\hat{\varphi}_2(r, 0) \neq 0.
\end{equation*}
\end{Hypothesis}
Assume Hypothesis~\ref{h:PtoE2-non-snake}. Then Hypothesis \ref{h:unique} means that the gradient of $g(\cdot,0,\cdot)$ does not vanish at any point within the set $\Gamma$.
\begin{Hypothesis}\label{h:hc-unique}
\begin{equation*}
\{\hat\varphi(r,0) - h^c(0,\hat{\lambda}_2(r,0)),  \, r \in S^1 \}\subsetneq S^1.
\end{equation*}
\end{Hypothesis}
The subtraction in the hypothesis is done in $S^1$. Since $I_\varphi$ is a
proper subset of $S^1$, this hypothesis is satisfied if
$|D_{\lambda_2}	h^c(0,\hat{\lambda}_2(r,0)) D_r\hat{\lambda}_2(r, 0)|\ll 1,$ $\forall \, r \in S^1$.
\begin{Theorem}\label{t:nonsnaking}
Assume Hypotheses \ref{h:general}-\ref{h:PtoE1}, \ref{h:EtoP} $(i)$ and Hypotheses~\ref{h:PtoE2-non-snake}-\ref{h:hc-unique}. Then, there is a sequence of mutually disjoint 
closed continuation curves 
$\Lambda_k:=\{(\lambda_{1,k}(r),\lambda_{2,k}(r)), r\in S^1\}$, $k\in\N$.
These curves tend towards a segment of the $\lambda_2$-axis in the sense of the Hausdorff metric.
More precisely,
for sufficiently large $k \in \N$, there exist mutually disjoint intervals ${\mathcal I}_k \subset \R$ and 
 smooth functions $\lambda_{i,k} \, : \, S^1 \rightarrow J_i, \ i=1,2$, and $L_k:S^1\to {\mathcal I}_k$
such that there is a homoclinic orbit to $E$ for $\lambda \in Q$ with flight time $L \in {\mathcal I}_k$ from $\Sigma^{in}$ to $\Sigma^{out}$, if and only if, there exists an $r \in S^1$ such that $\lambda=(\lambda_{1,k}(r), \, \lambda_{2,k}(r))$ and $L=L_k(r)$.
Moreover
\begin{align*}
 		\lambda_{1,k}(r) &=  O(e^{-\eta k})\cr
 		\lambda_{2,k}(r)&=  \hat{\lambda}_2(r,\lambda_{1,k}(r) ) + O(e^{-\eta k}) .
 \end{align*}
\end{Theorem}

\begin{figure}[ht]
\centering
\input{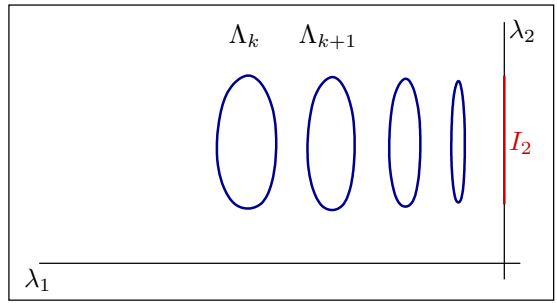}
 \caption{Visualisation of Theorem~\ref{t:nonsnaking}: The continuation curves $\Lambda_k$.}\label{f:nonsnaking_curve}
\end{figure}

\begin{prooof}[\,\,of Theorem~\ref{t:nonsnaking}]

We follow the lines of the proof of Theorem~\ref{t:snaking} up to equations  \eqref{e:coup_in_hc_rew}~--~\eqref{e:coup_out_rewr}, which we repeat here:
\begin{align}
 \varphi&=h^c(\lambda_1,\lambda_2),\label{e:coup_in_hc_rew-non-snake}
\\
a_u e^{-2\alpha^u(\lambda_1,\lambda_2)l}(1+O(e^{-\eta l}))&=h^u(\lambda_1,\lambda_2),\label{e:coup_in_hu_rew-non-snake}
\\
a_s e^{2\alpha^s(\lambda_1,\lambda_2)l}(1+O(e^{-\eta l}))&=g(\varphi+2l,\lambda_1,\lambda_2).\label{e:coup_out_rew-non-snake}
\end{align}
Again we look for solutions of these equations near the set $\Gamma$. For sufficiently small $\hat\lambda_1 \in J_1$
each point of a small tubular neighbourhood of ${\Gamma}\cap\{\lambda_1=\hat\lambda_1\}$ has a unique representation:
\begin{align*}
(\varphi + 2l)\!\!\!\mod 2\pi &= \hat{\varphi}(r,\lambda_1) + \tilde{\varphi} \cr
\lambda_2 &= \hat{\lambda}_2(r, \lambda_1) +\tilde{\lambda}_2,
\end{align*}
where $(\tilde{\varphi}, \tilde{\lambda}_2)$ is in the normal space
$\{(D_r\hat{\varphi}(r, \lambda_1),  \, D_r\hat{\lambda}_2(r, \lambda_1))
\}^\bot$ of $\Gamma$ at the point $(\hat{\varphi}(r,\lambda_1),
\hat{\lambda}_2(r, \lambda_1) )$. Hence, $\tilde{\varphi}$ and
$\tilde{\lambda}_2$ have to satisfy the additional equation:
\begin{equation}\label{e:coup_in_n_gam}
D_r\hat{\varphi}(r,\lambda_1) \, \tilde{\varphi} +  D_r\hat{\lambda}_2(r, \lambda_1) \, \tilde{\lambda}_2 =0.
\end{equation}
At first we solve equation \eqref{e:coup_out_rew-non-snake} together with equation \eqref{e:coup_in_n_gam}. Expanding $g$ w.r.t. $(\tilde{\varphi}, \tilde{\lambda}_2)$ we obtain
\begin{equation*}
g(\varphi +2l, \lambda_1, \lambda_2) = g(\hat{\varphi},\lambda_1,\hat{\lambda}_2)+ D_{\varphi}g(\hat{\varphi},\lambda_1,\hat{\lambda}_2) \tilde{\varphi} + D_{\lambda_2}g(\hat{\lambda}_2,\lambda_1,\hat{\lambda}_2) \tilde{\lambda}_2 + O(|(\tilde{\varphi} , \tilde{\lambda}_2)|^2).
\end{equation*}
Inserting this into \eqref{e:coup_out_rew-non-snake}
(and in the process exploiting \eqref{e:Gamma_rep-non-snake}) yields 
\begin{equation}\label{e:coup-1-non-snake}
D_{\varphi}	g(\hat{\varphi},\lambda_1,\hat{\lambda}_2) \tilde{\varphi} + D_{\lambda_2}	g(\hat{\varphi},\lambda_1,\hat{\lambda}_2)  \tilde{\lambda}_2
+ {O}(|(\tilde{\varphi} , \tilde{\lambda}_2)|^2)- a_s e^{2\alpha^s(\lambda_1,\hat{\lambda}_2(r, \lambda_1) +\tilde{\lambda}_2)l}(1+O(e^{-\eta l})) =0.
\end{equation}

Because of Hypothesis \ref{h:unique}, we may apply the implicit function theorem
to the system (\eqref{e:coup_in_n_gam},
\eqref{e:coup-1-non-snake}), and find a
unique solution $(\tilde{\varphi}^*,\tilde{\lambda}_2^*)(l,r,\lambda_1)$
for sufficiently small $\lambda_1$ and $|(\tilde{\varphi} ,
\tilde{\lambda}_2)|$ and sufficiently large $l$. This solution tends
to zero uniformly in $r$ and $\lambda_1$ as
$l$ tends to infinity.  Moreover,
$(\tilde{\varphi}^*,\tilde{\lambda}_2^*)(l,r,\lambda_1)$ is differentiable
w.r.t. $l$ and $D_l(\tilde{\varphi}^*,\tilde{\lambda}_2^*)(l,r,\lambda_1) =
O(e^{-\eta l})$. The latter can be seen by differentiating
\eqref{e:coup-1-non-snake} 
and taking Lemma  \ref{l:shilnikov} into consideration.

Next we consider \eqref{e:coup_in_hu_rew-non-snake}, where we
insert $(\tilde{\varphi}^*,\tilde{\lambda}_2^*)$.  Due to Hypothesis
\ref{h:EtoP}~(i), the resulting equation can be solved in the same way as
\eqref{e:coup_in_hu_rew} in Section \ref{s:snak_ana}.
We obtain the solution $\lambda_1=\lambda_1(l,r)$, and we find that both
$\lambda_1(l,r)$ and $ D_l \lambda_1(l,r)$ are of order $O(e^{-\eta l})$. The
estimates for the derivatives follow from Lemma~\ref{l:shilnikov}.

It remains to solve equation \eqref{e:coup_in_hc_rew-non-snake}, which can be written as
\begin{equation}\label{e:hc-non-snake}
(\hat{\varphi}(r,\lambda_1(l,r)) + \tilde{\varphi}(l,r) - 2l)\!\!\!\mod 2\pi = h^c(\lambda_1(l,r), \hat{\lambda}_2(r,\lambda_1(l,r)) + \tilde{\lambda}_2(l,r)),
\end{equation}
where $(\tilde{\varphi},\tilde{\lambda}_2)(l,r):=(\tilde{\varphi}^*,\tilde{\lambda}_2^*)(l,r,\lambda_1^*)$.
To solve \eqref{e:hc-non-snake}, we have to
overcome similar obstacles as when solving equation \eqref{e:coup_in_hc_rewr}. 
For that purpose we define 
\begin{equation}\label{e:l_0_def}
 2l_0(r):= \hat{\varphi}(r,0)  - h^c(0, \hat{\lambda}_2(r,0)).
\end{equation}
Note that $l_0(r): S^1\to \R$ is smooth.
Now we set
$2l= 2l_0(r) + 2\tilde{l} +2 k \pi$, for some $k \in \N$. 
Fixing $k \in \N$, we define
\begin{equation*}
 \lambda_{1,k}(\tilde{l},r):=\lambda_1(l_0(r) + \tilde{l} + k \pi,r),\quad
\lambda_{2,k}(\tilde{l},r):=\lambda_2(l_0(r) + \tilde{l} + k \pi,r),\quad
\tilde{\varphi}_k(\tilde{l},r):=\tilde{\varphi}(l_0(r) + \tilde{l} + k \pi,r).
\end{equation*}
Using this, we rewrite \eqref{e:hc-non-snake} as
\begin{equation}\label{e:hc2-non-snake}
 \tilde{\varphi}_k(\tilde{l},r)+ O(\lambda_{1,k})+ O(\tilde{\lambda}_{2,k})  =   2\tilde{l}.
\end{equation}
This equation can be solved by means of the contraction principle. For that we
note that the terms on the left-hand side together with their derivatives are of
order $O(e^{-\eta k })$.
Thus, for each fixed $r$ and sufficiently large $k$, equation
\eqref{e:hc2-non-snake} possesses a unique fixed point
$\tilde{l}^{*}_k(r)$. Moreover, $\tilde{l}^{*}_k: S^1\to\R$ is smooth, and
$\tilde{l}^{*}_k(r)$ is of order $O(e^{-\eta k })$.
All in all we obtain the unique solutions 
 \begin{align}
 		 l_k(r) &= {l}_0(r) + \tilde{l}_k(r) + k \pi \label{e:l_k_def},
\\
 		\lambda_{1,k}(r) &=\lambda_{1}({l}_k(r),r), \nonumber
\\
 		\lambda_{2,k}(r)&=  \hat{\lambda}_2(r,\lambda_{1,k}(r) ) + \tilde{\lambda}_2({l}_k(r),r).\nonumber
 \end{align}
Obviously, $\Lambda_k:=\{(\lambda_{1,k}(r),
\lambda_{2,k}(r)), r \in S^1\}$ are closed curves. Finally, with the
transformation 
\begin{equation}\label{e:L_transf_re_r}
 L_k(r)=l_k(r) + O(e^{-\eta l_k(r)}),
\end{equation}
cf. \eqref{e:L_transf_re}, we get
 \begin{align*}
 		\lambda_{1,k}(r) &=  \lambda_{1}(L_k(r) + O(e^{-\eta \, l_k(r)}),r) \, \ = \ \, O(e^{-\eta k}),\cr
 		\lambda_{2,k}(r)&=  \hat{\lambda}_2(r,\lambda_{1,k}(r) ) + \tilde{\lambda}_2(L_k(r) + O(e^{-\eta l_k(r)}),r)  \, \ = \ \, \hat{\lambda}_2(r,\lambda_{1,k}(r) ) + O(e^{-\eta k}) .
 \end{align*}

 Furthermore, we define the intervals ${\mathcal I}_k=[\underline{l}_k, \overline{l}_k]:= L_k(S^1)$. Due to Hypothesis~\ref{h:hc-unique} and \eqref{e:l_0_def}, the length of the interval
 ${l}_0(S^1)$ is less than $\pi$. Finally, since
 $\tilde{l}_k(r) = O(e^{-\eta k })$, it follows with \eqref{e:l_k_def} and \eqref{e:L_transf_re_r} that there is a $d>0$ such that for sufficiently large $k$
\begin{equation}\label{e:intersec_I_k}
 \underline{l}_{k+1}-\overline{l}_k>d.
\end{equation}
Hence, for sufficiently large $k$ and $\hat k$
\begin{equation*}
 {\mathcal I}_k\cap {\mathcal I}_{\hat k}=\emptyset .
\end{equation*}

It remains to show that the curves $\Lambda_k$ are mutually disjoint:
Assume that there exist $k, \hat{k} \in \N, \, \hat{k}>k $ and $r, \hat{r} \in S^1$ such that
\[(\lambda_1,\lambda_2):=(\lambda_{1,k}(r),\lambda_{2,k}(r))= (\lambda_{1,\hat{k}}(\hat{r}),\lambda_{2,\hat{k}}(\hat{r})). \]
Hence,
$h^u(\lambda_1,\lambda_2)=h^u(\lambda_{1,k}(r),\lambda_{2,k}(r))= h^u(\lambda_{1,\hat{k}}(\hat{r}),\lambda_{2,\hat{k}}(\hat{r}))$,
and from
\eqref{e:coup_in_hu_rew-non-snake} we deduce
\begin{equation}\label{e:l_k}
e^{-2\alpha^u(\lambda_1,\lambda_2)l_k}(1+O(e^{-\eta l_k})) =
e^{-2\alpha^u(\lambda_1,\lambda_2)l_{\hat{k}}}(1+O(e^{-\eta l_{\hat{k}}})),
\end{equation}
where we exploited the fact that $a_u$ depends only on $\lambda_1, \lambda_2$
and $\varphi$, together with equation \eqref{e:coup_in_hc_rew-non-snake}. Because of $l_{\hat{k}}(\hat r)- l_k(r)>d$, cf. \eqref{e:intersec_I_k}, we infer from \eqref{e:l_k} that
\begin{equation}\label{e:l_k2}
(1+O(e^{-\eta l_k})) =
e^{-2\alpha^u(\lambda_1,\lambda_2)[l_{\hat{k}}-l_k]}(1+ O(e^{-\eta l_{\hat{k}}})) < e^{-2\alpha^u(\lambda_1,\lambda_2)d}+e^{-2\alpha^u(\lambda_1,\lambda_2)d} O(e^{-\eta l_{\hat{k}}}).
\end{equation}
Taking the limit $k, \hat k \rightarrow \infty$, we see that the left-hand side of \eqref{e:l_k2} tends to $1$.  On the other hand the right-hand side of \eqref{e:l_k2} is close to $e^{-2\alpha^u(\lambda_1,\lambda_2)d}<1$. This yields a contradiction, and for this reason the curves $\Lambda_k$ and $\Lambda_{\hat{k}}$ cannot intersect. 
\end{prooof}

\section{Numerical verification of the hypotheses}\label{s:numver}

In this section we show numerically that the assumptions stated in
section~\ref{s:setup_res} are satisfied for our motivating
example~\eqref{e:laser}.  Recall that the setting of system~\eqref{e:laser} is
exactly opposite to the setting described in the hypotheses in
Section~\ref{s:setup_res}, i.e. $\dim W^s(E) = 1$, meaning that the orbit
$\gamma_{\scriptscriptstyle{\rm PtoE}}$ is of codimension one and the orbit
$\gamma_{\scriptscriptstyle{\rm EtoP}}$ is of codimension zero.
In order to adapt the following numerical computations to the
explanations given in the preceding sections we consider the vector field
$-F(x,y,\varphi,\nu_1,\nu_2)$. Further we define
\begin{equation*}
 u:=(x,y,\varphi),\quad \lambda := (\nu_1,\nu_2),\quad f(u,\lambda):=-F(u,\lambda).
\end{equation*}
The statements of Hypothesis~\ref{h:EtoP} are
trivially satisfied:  As shown in~\cite{KrauRie:08}, the EtoP connection 
which exists for parameter values on $c_b$, cf. Figure~\ref{f:laser},
is detected by constructing a numerical test function that is defined by a
signed distance of the endpoint of an orbit segment and the starting point of
another orbit segment.  Each zero of the test function then corresponds to a
real EtoP connecting orbit.  It turns out that the roots of the test function
are indeed regular, which means that Hypothesis~\ref{h:EtoP} is numerically
satisfied.  What is more, the way the computations are set up immediately
provide the Floquet multipliers of the periodic orbit, hence we can easily check
that Hypothesis~\ref{h:general} is also satisfied (see also below for the
computation of the Floquet multipliers).  In the following we focus on the
numerical verification of Hypotheses~\ref{h:PtoE2} and~\ref{h:PtoE3}, which are
not trivial to check.

In our further considerations we restrict ourselves to verify the shape of
$\Gamma$ rather than to present a numerical verification of all requirements
stated in Hypotheses~\ref{h:PtoE1}~--~\ref{h:PtoE3}. However, the computation of
the shape of $\Gamma$
requires advanced numerical techniques.  In the following section, we show how
to do this
by finding and continuing the codimension-$0$ PtoE heteroclinic connections
$\gamma_{\scriptscriptstyle{\rm PtoE}}$ for parameter values $(\nu_1,\nu_2)$ along the curve $c_b$
(between the two intersection points with $t_0$, cf.  Figure~\ref{f:laser}).
More precisely, we compute the intersection points of
$\gamma_{\scriptscriptstyle{\rm PtoE}}$ with
a torus $\mathcal{T}=P\times S^1$ surrounding $P$ and
transform the coordinates of these intersection points to the required format.
Note that this means that $\Gamma$ appears as a curve similar to
the sketch in Figure~\ref{f:graph_g}.

\subsubsection*{The implementation of the method}
\label{sec:implementation}

In the following, we explain the different continuation runs (i.e. numerical
solutions of boundary value problems with varying parameters) needed for the
computation of the curve $\Gamma$ as defined in~\eqref{e:Gamma_def} for
system~\eqref{e:laser}.  
For the actual computations, we
utilise the software package \textsc{AUTO}~\cite{DoeAUTO07}, which requires us
to formulate the vector fields in the time-rescaled form $\dot{u} = T
f(u,\lambda)$, $T\in[0,1]$.

\subsubsection*{Step 1: Finding $\gamma_{\scriptscriptstyle{\rm PtoE}}$}

Similar to the computation in~\cite{KrauRie:08}, it is possible to find the
codimension-$0$ PtoE connecting orbit $\gamma_{\scriptscriptstyle{\rm PtoE}}$ by setting up and solving
appropriate boundary value problems.

We start by constructing orbit segments $u^+ \subset W^s(E)$ and $u^- \subset
W^u(P)$ (both reasonably close to $\gamma_{\scriptscriptstyle{\rm PtoE}}$) such
that $u^+(0),u^-(1) \in
\Sigma$, where $\Sigma$ is a cross-section of $\gamma_{\scriptscriptstyle{\rm PtoE}}$ dividing the
phase-space such that $E$ and $P$ are separated.  Then we close the gap $u^+(0)
- u^-(1)$, which corresponds to finding a numerical representation of
$\gamma_{\scriptscriptstyle{\rm PtoE}}$.  Note that this methods allows us to find both
possible PtoE connections.

In the following, the boundary value problems that are used for the consecutive
continuation runs with \textsc{AUTO} are listed.  Note that the equilibrium
point $E$ as well as its eigenspaces are analytically known and hence are
omitted in this listing.  

In order to compute the orbit segments, we need to continue several objects
simultaneously.  This can be achieved by extending the system by the additional
objects in consecutive continuation runs.

We assume that we have prepared the system such that $\lambda = (\nu_1,\nu_2)\in
c_b$, half way between the two intersection points with $t_b$.  Moreover, we
assume that the numerical representation of $P$ and the numerical representation
of the unstable Floquet bundle along $P$ and the value of the unstable Floquet
multiplier, which we here denote by $\mu$, is
known.  For more information about the necessary steps to achieve this,
we refer to \cite{KrauRie:08}.

The first object we need to continue is the numerical representation $u_P$ of the periodic orbit $P$.  The following
standard boundary value problem (BVP) for periodic orbits is used:
\begin{align}
  \label{eq:ode_p} \dot{u}_P &= T f(u_P, \lambda),\\
  \label{eq:bc_p} 0 & = u_P(0) - u_P(1),\\
  \label{eq:int_p} 0 & = \int_0^1 \langle \dot{\tilde{u}}_P(\tau), u_P(\tau)
  \rangle d\tau.
\end{align}
Note that during the continuation run, $\tilde{u}_P$ denotes a solution of $u_P$
from the previous computation step.

In addition to $u_P$, the Floquet bundle $u_F$ is continued using:
\begin{align}
  \label{eq:ode_f} \dot{u}_F &=  T D_u f(u_P(t), \lambda) u_F(t) + \ln|\mu|
  u_F(t),\\
  \label{eq:bc1_f} 0 &= \mbox{sgn}(\mu) u_F(0) - u_F(1),\\
  \label{eq:bc2_f} 0 &= \langle u_F(0), u_F(0) \rangle - 1.
\end{align}
Note that we assume that the Floquet bundle is normalised such that $\lVert
u_F(0) \rVert = 1$.  For more information on this BVP, see~\cite{Doedel2008}.

The orbit segment $u^- \subset W^u(P)$ is continued using the following
equations:
\begin{align}
  \label{eq:ode_minus} \dot{u}^- &= T^- f(u^-,\lambda),\\
  \label{eq:bc1_minus} 0 &= u_P(0) + \delta u_F(0) - u^-(0),\\
  \label{eq:bc2_minus} 0 &= \langle u^-(1) - \sigma, n_\Sigma \rangle - \eta^-.
\end{align}
Here, $\sigma$ is an arbitrary point in $\Sigma$ and $n_\Sigma$ is a normal of the cross-section
$\Sigma$.  The parameter $\delta$ is initialised with $-10^{-4}$, the solution
$u^-$ is initialised with a constant value of $u_P(0)+\delta u_F(0)$.

The first continuation run uses the ODEs~\eqref{eq:ode_p}, \eqref{eq:ode_f} and
\eqref{eq:ode_minus}.  The boundary conditions are \eqref{eq:bc_p},
\eqref{eq:bc1_f}, \eqref{eq:bc2_f}, \eqref{eq:bc1_minus}, \eqref{eq:bc2_minus},
the integral condition is \eqref{eq:int_p}.  The primary continuation parameter
is $T^-$ (initialised with $0$), the remaining continuation parameters are
$\mu$, $T$ and $\eta^-$.  The value of $\eta^-$ is considered a signed distance
of the endpoint $u^-(1)$ to $\Sigma$, hence if $\eta^- = 0$ the orbit segment
$u^-$ ends in $\Sigma$ and the continuation run stops.

The second orbit segment $u^+ \subset W^s(E)$ is defined by the following
equations:
\begin{align}
  \label{eq:ode_plus} \dot{u}^+ &= T^+ f(u^+,\lambda),\\
  \label{eq:bc1_plus} 0 &= E + \varepsilon v_s - u^+(1),\\
  \label{eq:bc2_plus} 0 &= \langle u^+(0) - \sigma, n_\Sigma \rangle - \eta^+.
\end{align}
Here, $v_s$ denotes a vector within the two-dimensional stable eigenspace at
$E$, which is analytically known.  The parameter $\varepsilon$ is initialised with
$10^{-6}$, the solution $u^+$ is initialised with a constant value of
$E+\varepsilon v_s$.  Note that we keep the boundary
condition~\eqref{eq:bc1_plus} throughout all following continuation runs, which
is only valid because $E$ is a saddle-focus.

The second continuation run uses the ODEs~\eqref{eq:ode_p}, \eqref{eq:ode_f},
\eqref{eq:ode_minus} and~\eqref{eq:ode_plus}.  The boundary conditions are
\eqref{eq:bc_p}, \eqref{eq:bc1_f}, \eqref{eq:bc2_f}, \eqref{eq:bc1_minus},
\eqref{eq:bc2_minus}, \eqref{eq:bc1_plus}, \eqref{eq:bc2_plus}, the integral
condition is \eqref{eq:int_p}.  The primary continuation parameter is $T^+$
(initialised with $0$), the remaining continuation parameters are $\mu$, $T$,
$T^-$ and $\eta^+$.  The value of $\eta^+$ is considered as a signed distance of
the starting point $u^+(0)$ to $\Sigma$, hence if $\eta^+ = 0$ the orbit segment
$u^+$ starts in $\Sigma$ and the continuation run stops. 

After these two continuation runs, the orbit segment $u^-$ ends in $\Sigma$ and
$u^+$ starts in $\Sigma$.  We define $z := (u^+(0) - u^-(1))/\lVert u^+(0) -
u^-(1)\rVert$ and initialise the new parameter $\eta$ with the value $\lVert
u^+(0) - u^-(1)\rVert$.

Then we replace the boundary condition~\eqref{eq:bc2_plus} by
\begin{align}\label{eq:bc_gap} 
0 &= u^+(0) - u^-(1) - \eta z,
\end{align}
which means that we force the difference $u^+(0) - u^-(1)$ to be in the linear
subspace defined by the vector $z$, while the parameter $\eta$ measures the gap
between these two points.

In order to close the gap, we perform a third continuation run using the
ODEs~\eqref{eq:ode_p}, \eqref{eq:ode_f}, \eqref{eq:ode_minus}
and~\eqref{eq:ode_plus}.  The boundary conditions are~\eqref{eq:bc_p},
\eqref{eq:bc1_f}, \eqref{eq:bc2_f}, \eqref{eq:bc1_minus}, \eqref{eq:bc2_minus},
\eqref{eq:bc1_plus}, \eqref{eq:bc_gap}, the integral condition is
\eqref{eq:int_p}.  The primary continuation parameter is $\eta$, the remaining
continuation parameters are $\mu$, $T$, $T^-$ and $T^+$.  The value of $\eta$ is
considered as a signed distance between the two points $u^+(0)$ and $u^-(1)$, hence
finding a root of $\eta$ corresponds to finding a numerical representation of
$\gamma_{\scriptscriptstyle{\rm PtoE}}$.  Note that there are two incarnations
of $\gamma_{\scriptscriptstyle{\rm PtoE}}$, which correspond to two different
roots of $\eta$.  This concludes the first step of our method.

\subsubsection*{Step 2: Computing the intersection point $\gamma_{\scriptscriptstyle{\rm PtoE}} \cap
\mathcal{T}$}

In this step we compute the intersection point of
$\gamma_{\scriptscriptstyle{\rm PtoE}}$ with the torus $\mathcal{T}$.  We use
the resulting orbit segments from step 1, but we no longer force the points
$u^-(1)$ and $u^+(0)$ to be in $\Sigma$, instead we let $u^-(1)$ vary along
$\gamma_{\scriptscriptstyle{\rm PtoE}}$.

For technical reasons, we need to include a second copy of $P$ into the
continuation using the equations:
\begin{align}
  \label{eq:ode_hatp} \dot{\hat{u}}_P &= \hat{T} f(\hat{u}_P,\lambda),\\
  \label{eq:bc_hatp} 0 &= \hat{u}_P(0) - \hat{u}_P(1),\\
  \label{eq:int_hatp} 0 & = \int_0^1 \langle \dot{\tilde{\hat{u}}}_P(\tau),
  \hat{u}_P(\tau) \rangle d\tau.
\end{align}
The solution $\hat{u}$ is initialised with the solution $u_P$, and the parameter
$\hat{T}$ is initialised with the parameter $T$.  Recall that the parameters $T$,
$\hat{T}$, $T^-$ and $T^+$ are necessary to transform the vector field such that
the time interval on which the orbits are computed is always $[0,1]$.

The following boundary condition is used to measure the distance of $u^-(1)$ to the Poincare section of $P$ at $\hat{u}_P(0)$:
\begin{align}
  \label{eq:bc_section} 0 &= \langle \hat{u}_P(0) - u^-(1),
  f(\hat{u}_P(0),\lambda) \rangle - \eta_1.
\end{align}

The fourth continuation run uses the ODEs~\eqref{eq:ode_p}, \eqref{eq:ode_f},
\eqref{eq:ode_minus}, \eqref{eq:ode_plus} and~\eqref{eq:ode_hatp}.  The boundary
conditions are~\eqref{eq:bc_p}, \eqref{eq:bc1_f}, \eqref{eq:bc2_f},
\eqref{eq:bc1_minus}, \eqref{eq:bc1_plus},
\eqref{eq:bc_gap}, \eqref{eq:bc_hatp}, \eqref{eq:bc_section}, the integral
conditions are~\eqref{eq:int_p} and~\eqref{eq:int_hatp}.  The primary
continuation parameter is $\eta_1$, the remaining continuation parameters are
$\mu$, $T$, $\hat{T}$, $T^+$, $T^-$, $\delta$, $\varepsilon$.  Any root of $\eta_1$
means that the point $u^-(1) = u^+(0)$ is in the Poincar{\'e} section defined by
$\hat{u}_P(0)$ and having the normal $f(\hat{u}_P(0),\lambda)$.  Note that in
general there may be several roots of $\eta_1$.  
We choose the solution for which $u^-(1)$ is in the unstable fibre of $\hat{u}_P(0)$, i.e. $u^-(1) \in W^{uu}(\hat{u}_P(0),\lambda)$.  In practise, we achieve this by using the solution for which the distance $\|u^-(1)-\hat{u}_P(0)\|$ is smallest.

In the final continuation run for the detection of the intersection point, we
measure the distance of $u^-(1)$ to $\hat{u}_P(0)$ using the boundary condition
\begin{align}
  \label{eq:bc_radius} 0 &= \lVert \hat{u}_P(0) - u^-(1) \rVert - \eta_2.
\end{align}

The fifth continuation run uses the ODEs~\eqref{eq:ode_p}, \eqref{eq:ode_f},
\eqref{eq:ode_minus}, \eqref{eq:ode_plus} and~\eqref{eq:ode_hatp}.  The boundary
conditions are~\eqref{eq:bc_p}, \eqref{eq:bc1_f}, \eqref{eq:bc2_f},
\eqref{eq:bc1_minus}, \eqref{eq:bc1_plus}, \eqref{eq:bc_gap},
\eqref{eq:bc_hatp}, \eqref{eq:bc_section}, \eqref{eq:bc_radius}, the integral
condition is~\eqref{eq:int_p}.  Note that we omit the integral phase condition
of $\hat{u}_P$.  The primary continuation parameter is $\eta_2$, the remaining
continuation parameters are $\mu$, $T$, $\hat{T}$, $T^+$, $T^-$, $\delta$,
$\varepsilon$.  When the value of $\eta_2$ reaches the desired radius of
$\mathcal{T}$, the continuation run stops and the intersection point is found.
The actual computations are performed for the radius $0.1$.

\subsubsection*{Step 3: Continuation of $\gamma_{\scriptscriptstyle{\rm PtoE}}$ along $c_b$}

Instead of using the (so far fixed) parameters $\nu_1$ and $\nu_2$ directly, we
define a smooth function $\lambda: [0,1] \to \R^2$, $s \mapsto \lambda(s)$, such
that $\lambda(0) = c_b \cap t_b$ (the lower intersection) and $\lambda(1) = c_b
\cap t_b$ (the upper intersection) and $\lambda([0,1]) \subset c_b$.  Using this
definition of $\lambda$, we can continue the system consisting
of~\eqref{eq:ode_p}, \eqref{eq:ode_f}, \eqref{eq:ode_minus}, \eqref{eq:ode_plus}
and~\eqref{eq:ode_hatp}, with boundary conditions~\eqref{eq:bc_p},
\eqref{eq:bc1_f}, \eqref{eq:bc2_f}, \eqref{eq:bc1_minus}, \eqref{eq:bc1_plus},
\eqref{eq:bc_gap}, \eqref{eq:bc_hatp}, \eqref{eq:bc_section},
\eqref{eq:bc_radius} and with integral condition~\eqref{eq:int_p}.
The primary continuation parameter is $s$ (initialised to $0.5$), the remaining
continuation parameters are $\mu$, $T$, $\hat{T}$, $T^+$, $T^-$, $\delta$,
$\varepsilon$.  This final continuation run is performed for increasing $s$
until a limit point for $s=1$ is reached, and then for decreasing $s$ until a
limit point for $s=0$ is reached.  These two limit points correspond to the two
intersection points of $c_b$ with the two branches of $t_b$.

\begin{figure}
  \begin{center}
    \includegraphics{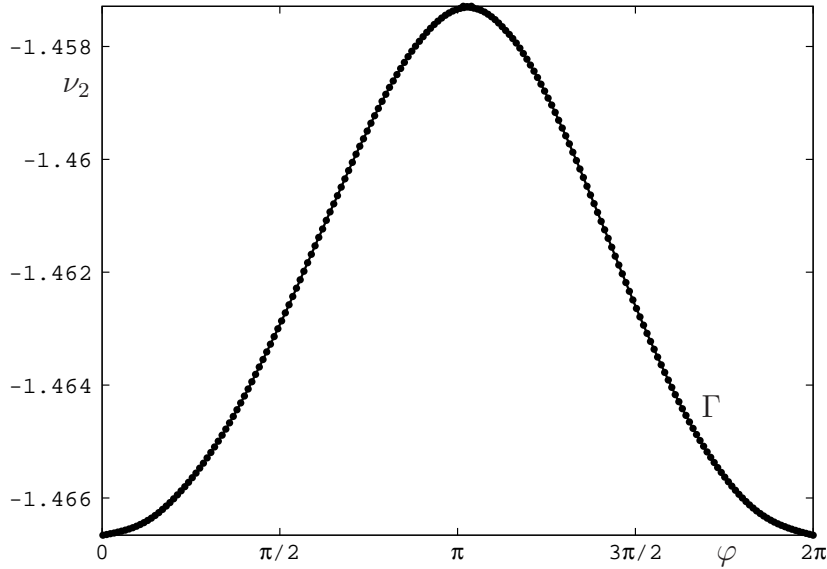}
  \end{center}
  \caption{The computed curve $\Gamma$ for system~\eqref{e:laser}.  The curve is
  shown as an $\varphi$-vs-$\nu_2$ plot, where the angle $\varphi$ is calculated
  as the local coordinate on the torus surface $\mathcal{T}$ and corresponds to
  the angle $\varphi$ used in Hypothesis~\ref{h:PtoE2}.  The curve $\Gamma$ is
  interpolated between the computed points (the marked points on $\Gamma$) using
  cubic splines and it is shifted such that the minimum value is at $\varphi =
  0$.\label{fig:numgamma}}
\end{figure}

This concludes our method of computing $\Gamma$ for system~\eqref{e:laser}, the
result after transforming the coordinates of $u^-(1) = u^+(0)$ to the local
coordinates of the torus is shown in Figure~\ref{fig:numgamma}.  Note that the
shape of $\Gamma$ verifies \ref{h:PtoE2}(ii) and \ref{h:PtoE3}, but it does not
yet cover \ref{h:PtoE1}(i).  In order to verify \ref{h:PtoE2}(i), the curves
corresponding to $\mbox{graph}\; g$ for different values of $\varphi$ (cf.
Figure~\ref{f:graph_hu}) needs to be computed.  The computations are very
similar to the computations presented above and are therefore omitted here.

\section{Discussion and conclusion}
In this paper we have considered the behaviour of one-homoclinic orbits near EtoP cycles. We have rigorously shown how the behaviour of the primary cycle determines the continuation behaviour of the homoclinic orbits.
Our analysis is restricted to $\R^3$.
In higher dimensionsal state space the dimensions of $\Sigma^{in}$ and $\Sigma^{out}$ increase accordingly and therefore also the number of bifurcation equations. 

In our considerations we distinguished the cases that the periodic orbit $P$ has positive or negative Floquet multipliers, respectively. 
For positive Floquet multipliers we have discussed two different scenarios. First we have verified homoclinic snaking as it was previously observed numerically in our motivating example \eqref{e:laser}. Further we have described a nonsnaking scenario which has to our knowledge not yet been observed in systems in $\R^3$. In systems in $\R^4$ however this effect was observed numerically. Though it is not clear whether or not in these examples this effect is due to the behaviour of the EtoP cycle as assumed in the present paper.
For negative Floquet multipliers we confined to study the corresponding snaking scenario.

For the detection of multi-around homoclinic orbit to $E$, these are orbits that follow the primary EtoP cycle several times before returning to $E$, also couplings near the equilibrium $E$ have to be considered. At least if these investigations are combined with higher dimensional systems a Lin's method approach seems to be appropriate. However, numerical results in \cite{Rie10} give rise to the hope that those homoclinic orbit exist in our motivating example.

In Section~\ref{s:numver} we verified numerically that
system~\eqref{e:laser} satisfies the hypotheses which we 
used in our analysis. In particular we used a novel approach 
based on numerical continuation techniques to compute 
$\Gamma$. The computational results are very satisfying, the 
computed shape of $\Gamma$ looks exactly like expected. This 
shows that Hypothesis~\ref{h:PtoE2}(ii) is true for our 
motivating system.

\paragraph{Acknowledgments.} M.V. was supported by the German National Academic Foundation (Studienstifting des deutschen Volkes). The authors are grateful to V. Kirk for pointing out related results.

\end{document}